
\def\LaTeX{%
  \let\Begin\begin
  \let\End\end
  \let\salta\relax
  \let\finqui\relax
  \let\futuro\relax}

\def\UK{\def\our{our}\let\sz s}
\def\USA{\def\our{or}\let\sz z}

\UK



\LaTeX

\USA


\salta

\documentclass[twoside,12pt]{article}
\setlength{\textheight}{24cm}
\setlength{\textwidth}{16cm}
\setlength{\oddsidemargin}{2mm}
\setlength{\evensidemargin}{2mm}
\setlength{\topmargin}{-15mm}
\parskip2mm


\usepackage[usenames,dvipsnames]{color}
\usepackage{amsmath}
\usepackage{amsthm}
\usepackage{amssymb}
\usepackage[mathcal]{euscript}
\usepackage{hyperref}
%
%


\definecolor{viola}{rgb}{0.3,0,0.7}
\definecolor{ciclamino}{rgb}{0.5,0,0.5}
\definecolor{rosso}{rgb}{0.85,0,0}

\def\anold #1{{\color{rosso}#1}}
\def\an #1{{\color{rosso}#1}}
\def\gianni #1{{\color{red}#1}}
\def\betti #1{{\color{blue}#1}}

\def\rev #1{{\color{rosso}#1}}

\def\betti #1{#1}
\def\juerg #1{#1}

\def\anold #1{#1}
\def\an #1{#1}
\def\rev #1{#1}




\bibliographystyle{plain}


%

\finqui

\def\Beq{\Begin{equation}}
\def\Eeq{\End{equation}}
\def\Bsist{\Begin{eqnarray}}
\def\Esist{\End{eqnarray}}

\def\Bthm{\Begin{theorem}}
\def\Ethm{\End{theorem}}
\def\Blem{\Begin{lemma}}
\def\Elem{\End{lemma}}

\def\Brem{\Begin{remark}\rm}
\def\Erem{\End{remark}}

\def\Bdim{\Begin{proof}}
\def\Edim{\End{proof}}
\def\Bcenter{\Begin{center}}
\def\Ecenter{\End{center}}
\let\non\nonumber




\def\step #1 \par{\medskip\noindent{\bf #1.}\quad}


\def\Lip{Lip\-schitz}
\def\Holder{H\"older}
\def\frechet{Fr\'echet}
\def\aand{\quad\hbox{and}\quad}

\def\lhs{left-hand side}
\def\rhs{right-hand side}
\def\sfw{straightforward}


\def\nbh{neighb\our hood}


\def\multibold #1{\def\arg{#1}%
  \ifx\arg\pto \let\next\relax
  \else
  \def\next{\expandafter
    \def\csname #1#1#1\endcsname{{\boldsymbol #1}}%
    \multibold}%
  \fi \next}

\def\pto{.}

\def\multical #1{\def\arg{#1}%
  \ifx\arg\pto \let\next\relax
  \else
  \def\next{\expandafter
    \def\csname cal#1\endcsname{{\cal #1}}%
    \multical}%
  \fi \next}

\def\multigrass #1{\def\arg{#1}%
  \ifx\arg\pto \let\next\relax
  \else
  \def\next{\expandafter
    \def\csname gr#1\endcsname{{\mathbb #1}}%
    \multigrass}%
  \fi \next}


\def\multimathop #1 {\def\arg{#1}%
  \ifx\arg\pto \let\next\relax
  \else
  \def\next{\expandafter
    \def\csname #1\endcsname{\mathop{\rm #1}\nolimits}%
    \multimathop}%
  \fi \next}

\multibold
qwertyuiopasdfghjklzxcvbnmQWERTYUIOPASDFGHJKLZXCVBNM.

\multical
QWERTYUIOPASDFGHJKLZXCVBNM.

\multigrass
QWERTYUIOPASDFGHJKLZXCVBNM.

\multimathop
diag dist div dom mean meas sign supp .

\def\Span{\mathop{\rm span}\nolimits}


\def\accorpa #1#2{\eqref{#1}--\eqref{#2}}
\def\Accorpa #1#2 #3 {\gdef #1{\eqref{#2}--\eqref{#3}}%
  \wlog{}\wlog{\string #1 -> #2 - #3}\wlog{}}


\def\separa{\noalign{\allowbreak}}

\def\neto{\mathrel{{\scriptscriptstyle\nearrow}}}
\def\seto{\mathrel{{\scriptscriptstyle\searrow}}}

\def\graffe #1{\mathopen\{#1\mathclose\}}

\def\<#1>{\mathopen\langle #1\mathclose\rangle}
\def\norma #1{\mathopen \| #1\mathclose \|}

\def\[#1]{\mathopen\langle\!\langle #1\mathclose\rangle\!\rangle}

\def\iot {\int_0^t}
\def\ioT {\int_0^T}
\def\intQt{\int_{Q_t}}
\def\intQ{\int_Q}
\def\iO{\int_\Omega}
\def\iG{\int_\Gamma}
\def\intS{\int_\Sigma}
\def\intSt{\int_{\Sigma_t}}

\def\dt{\partial_t}
\def\dn{\partial_{\boldsymbol\nu}}

\def\cpto{\,\cdot\,}

\def\checkmmode #1{\relax\ifmmode\hbox{#1}\else{#1}\fi}
\def\aeO{\checkmmode{a.e.\ in~$\Omega$}}
\def\aeQ{\checkmmode{a.e.\ in~$Q$}}
\def\aeG{\checkmmode{a.e.\ on~$\Gamma$}}
\def\aeS{\checkmmode{a.e.\ on~$\Sigma$}}
\def\aet{\checkmmode{a.e.\ in~$(0,T)$}}

\def\aat{\checkmmode{for a.a.~$t\in(0,T)$}}


\def\erre{{\mathbb{R}}}




\def\genspazio #1#2#3#4#5{#1^{#2}(#5,#4;#3)}
\def\spazio #1#2#3{\genspazio {#1}{#2}{#3}T0}

\def\L {\spazio L}
\def\H {\spazio H}
\def\W {\spazio W}

\def\C #1#2{C^{#1}([0,T];#2)}


\def\Lx #1{L^{#1}(\Omega)}
\def\Hx #1{H^{#1}(\Omega)}

\def\LxG #1{L^{#1}(\Gamma)}
\def\HxG #1{H^{#1}(\Gamma)}

\def\LQ #1{L^{#1}(Q)}
\def\LS #1{L^{#1}(\Sigma)}

\def\Luno{\Lx 1}
\def\Ldue{\Lx 2}
\def\Linfty{\Lx\infty}

\def\Huno{\Hx 1}
\def\Hdue{\Hx 2}

\def\HunoG{\HxG 1}
\def\HdueG{\HxG 2}

\def\LunoG{\LxG 1}
\def\LdueG{\LxG 2}


\let\badphi\phi
\let\theta\vartheta
\let\eps\varepsilon
\let\phi\varphi

\let\TeXchi\chi                         
\newbox\chibox
\setbox0 \hbox{\mathsurround0pt $\TeXchi$}
\setbox\chibox \hbox{\raise\dp0 \box 0 }
\def\chi{\copy\chibox}


\def\QED{\hfill $\square$}


\let\emb\hookrightarrow
\def\CO{C_\Omega}

\def\suG{_{|\Gamma}}

\def\misO{|\Omega|}
\def\misG{|\Gamma|}

\def\VG{V_\Gamma}
\def\HG{H_\Gamma}
\def\WG{W_\Gamma}
\def\nablaG{\nabla_\Gamma}
\def\DeltaG{\Delta_\Gamma}
\def\muG{\mu_\Gamma}
\def\phiG{\phi_\Gamma}
\def\xiG{\xi^\Gamma}
\def\FG{F_\Gamma}

\def\pG{p_\Gamma}
\def\qG{q_\Gamma}
\def\vG{v_\Gamma}
\def\wG{w_\Gamma}
\def\PG{P_\Gamma}
\def\QG{Q_\Gamma}

\def\fG{f^\Gamma}
\def\gG{g^\Gamma}
\def\hG{h^\Gamma}
\def\uG{u^\Gamma}
\def\zG{z^\Gamma}

\def\tarQ{\badphi^Q}
\def\tarS{\badphi^\Sigma}
\def\tarO{\badphi^\Omega}
\def\tarG{\badphi^\Gamma}

\def\psiG{\psi_\Gamma}
\def\etaG{\eta_\Gamma}
\def\rhoG{\rho_\Gamma}
\def\thetaG{\theta_\Gamma}

\def\phiz{\phi_0}
\def\phiGz{{\phiz}\suG}
\def\mz{m_0}
\def\rz{r_0}

\def\mueps{\mu^\eps}
\def\phieps{\phi^\eps}
\def\muGeps{\mu_\Gamma^\eps}
\def\phiGeps{\phi_\Gamma^\eps}
\def\aeps{\alpha_\eps}
\def\peps{p^\eps}
\def\qeps{q^\eps}
\def\pGeps{\pG^\eps}
\def\qGeps{\qG^\eps}

\def\un{u_n}
\def\uGn{u^\Gamma_n}
\def\mun{\mu^n}
\def\muGn{\mu_\Gamma^n}
\def\phin{\phi^n}
\def\phiGn{\phi_\Gamma^n}

\def\muh{\hat\mu}
\def\muGh{\hat\mu_\Gamma}
\def\phih{\hat\phi}
\def\phiGh{\hat\phi_\Gamma}
\def\chiG{\chi^\Gamma}

\def\ustar{u_*}
\def\uGstar{u^\Gamma_*}
\def\mustar{\mu^*}
\def\muGstar{\muG^*}
\def\phistar{\phi^*}
\def\phiGstar{\phiG^*}

\def\soluz{(\mu,\muG,\phi,\phiG,\xi,\xiG)}
\def\soluzeps{(\mueps,\muGeps,\phieps,\phiGeps)}
\def\soluzl{(\eta,\etaG,\psi,\psiG)}
\def\soluza{(p,\pG,q,\qG)}
\def\soluzaeps{(\peps,\pGeps,\qeps,\qGeps)}
\def\sol{(\mu,\muG,\phi,\phiG)}
\def\solstar{(\mustar,\muGstar,\phistar,\phiGstar)}
\def\solh{(\muh,\muGh,\phih,\phiGh)}

\let\lam\lambda
\def\lamG{\lam^\Gamma}
\let\Lam\Lambda
\def\LamG{\Lam^\Gamma}
\def\Beta{\hat\beta}
\def\BetaG{\Beta_\Gamma}
\def\betaG{\beta_\Gamma}
\def\betaeps{\beta_\eps}
\def\betaGeps{\beta_{\Gamma\!,\,\eps}}
\def\Betaeps{\hat\beta_\eps}
\def\BetaGeps{\hat\beta_{\Gamma\!,\,\eps}}
\def\betaz{\beta^\circ}
\def\betaGz{\betaG^\circ}
\def\Pi{\hat\pi}
\def\PiG{\Pi_\Gamma}
\def\piG{\pi_\Gamma}

\def\phimin{\phi_*}
\def\phimax{\phi^*}
\def\umin{u_{\an{\min}}}
\def\umax{u_{\an{\max}}}
\def\uGmin{u^\Gamma_{\min}}
\def\uGmax{u^\Gamma_{\max}}
\def\Uad{\calU_{ad}}
\def\UR{\calU_R}

\def\calVz{\calV_0}
\def\calVp{\calV^{\,*}}

\def\normaH #1{\norma{#1}_H}

\def\normaHG #1{\norma{#1}_{\HG}}

\def\normaHH #1{\norma{#1}_{\calH}}
\def\normaVV #1{\norma{#1}_{\calV}}
\def\normaWW #1{\norma{#1}_{\calW}}

\let\hat\widehat

\Begin{document}


%
\title{Well-posedness and optimal control\\[0.3cm] 
  for a viscous Cahn--Hilliard--Oono system\\[0.3cm] 
  with dynamic boundary conditions}
\author{}
\date{}
\maketitle
\Bcenter
\vskip-1cm
{\large\sc Gianni Gilardi$^{(1)(2)}$}\\
{\normalsize e-mail: {\tt gianni.gilardi@unipv.it}}\\[.25cm]
{\large\sc Elisabetta Rocca$^{(1)(2)}$}\\
{\normalsize e-mail: {\tt elisabetta.rocca@unipv.it}}\\[.25cm]
{\large\sc Andrea Signori$^{(3)}$}\\
{\normalsize e-mail: {\tt andrea.signori@polimi.it}}\\[0.25cm]
$^{(1)}$
{\small Dipartimento di Matematica ``F. Casorati'', Universit\`a di Pavia}\\
{\small via Ferrata 5, 27100 Pavia, Italy}\\[.2cm]
$^{(2)}$
{\small Research Associate at the IMATI -- C.N.R. Pavia}\\[.2cm]
$^{(3)}$
{\small Dipartimento di Matematica, Politecnico di Milano}\\
{\small via Bonardi 9, 20133 Milano, Italy}\\[.2cm]
\Ecenter
\begin{center}
\emph{Dedicated to our dear friend Pierluigi Colli\\
on the occasion of his 65th birthday}
\end{center}

\Begin{abstract}
\noindent
\betti{In this paper we consider a nonlinear system of PDEs coupling the viscous Cahn--Hilliard--Oono equation with dynamic boundary conditions \anold{enjoying a similar structure on the boundary}. After proving well-posedness of the {corresponding} initial \anold{boundary} value problem, we study \anold{an associated} optimal control problem related to a \an{tracking-type} cost functional, proving \anold{first-order} necessary conditions of optimality. \anold{The controls enter \an{the} system in the form of a distributed and a boundary source.} We can account for general potentials in the bulk and in the boundary part \an{under the common assumption} that the boundary potential is dominant with respect to the \anold{bulk one}.
For example, \anold{the} regular quartic potential as well as \anold{the logarithmic potential} can be considered in our analysis.}
\vskip3mm
\noindent {\anold{\bf Keywords:}}
\anold{
Cahn--Hilliard--Oono equation, dynamic boundary conditions, bulk-surface mass dynamics, well-posedness, optimal control.
}
\vskip3mm
\noindent {\bf AMS (MOS) Subject Classification:} 
\anold{
{
		35K55, 
        35K61, 
		49J20, 
		49J50, 
		49K20. 
		}
}
\End{abstract}
\salta
\pagestyle{myheadings}
\newcommand\testopari{\sc Gilardi \ --- \ Rocca \ --- \ Signori}
\newcommand\testodispari{{\sc \anold{ bulk-surface Cahn--Hilliard--Oono dynamics}}}
\markboth{\testodispari}{\testopari}
\finqui
%

\section{Introduction}
\label{Intro}
\setcounter{equation}{0}

\betti{In this paper, we study the following nonlinear PDE system, referred to in the literature as the viscous Cahn--Hilliard--Oono (CHO) system (cf.~\cite{oonopuri87, oonopuri88I, oonopuri88II}), which is of great interest in the study of pattern formations in 
phase-separating materials} 
\Beq
  \dt\phi - \Delta\mu = \gamma(u-\phi)
  \aand
  \tau \dt\phi - \Delta\phi + F'(\phi) = \mu
  \quad \hbox{in $Q:=\Omega\times(0,T)$}.
  \label{Isystem}
\Eeq
\betti{Here the unknowns are the order parameter~$\phi$ and the chemical potential~$\mu$, while} 
 $\Omega$ is the \anold{spatial} domain where the evolution takes place and $T$ is a \an{prescribed} final time.
In the above equations, 
$\gamma$~is a positive constant,
$u$~is a \anold{source/sink term that later on will play the role of (distributed) control},
$\tau$~is a \anold{positive} viscosity coefficient,
and $F'$ is the derivative of a double-well potential~$F$.

Typical and important examples \anold{for this latter}
are the so-called classical  {\em regular potential} and the {\em logarithmic} double-well potential.
They are given\an{, in the order,} by
\begin{align}
  & F_{\an{ \rm reg}}(r) := \frac 14 \, (r^2-1)^2 \,,
  \quad r \in \erre, 
  \label{regpot}
  \\
  & F_{\an{\rm log}}(r) := ((1+r)\ln (1+r)+(1-r)\ln (1-r)) - \anold{c_1} r^2 \,,
  \quad r \in (-1,1),
  \label{logpot}
\end{align}
where $\anold{c_1}>1$ \juerg{is such} that $F_{\an{\rm log}}$ is nonconvex
{(then, $F_{log}$ is extended outside the physical interval $(-1,1)$ in the usual way,
i.e., by continuity at $\pm1$ and as $+\infty$ elsewhere,
in order to preserve semicontinuity).}
Another example is the following {\em double obstacle potential\/}, where $\anold{c_2}>0$,
\Beq
  F_{\an{\rm 2obs}}(r) := - \anold{c_2}r^2 
  \quad \hbox{if $|r|\leq1$}
  \aand
  F_{\an{\rm 2obs}}(r) := +\infty
  \quad \hbox{if $|r|>1$}.
  \label{obspot}
\Eeq
In cases like \eqref{obspot}, one has to split $F$ into a nondifferentiable convex part 
(the~indicator function of $[-1,1]$ in the present example) and a smooth perturbation \an{(typically quadratic)}.
Accordingly, one has to replace the derivative of the convex part \an{of the potential}
by \an{its} subdifferential and interpret the second identity in \eqref{Isystem} as a differential inclusion.

\betti{
System \eqref{Isystem} can be seen as a Cahn--Hilliard (CH) equation with reaction and
turns out to be useful in several applications
such as biological models~\cite{KS}, inpainting algorithms~\cite{BEG}, and
polymers~\cite{BO}.
The main feature that makes the CHO equation different from the standard CH one is the fact that in \eqref{Isystem} the conservation of mass is not guaranteed even in presence of ``standard'' Neumann homogeneous boundary conditions \anold{for the chemical potential $\mu$}.
The reaction term in \eqref{Isystem} takes into account long-range interactions and, of course, more general terms could be considered, in particular\an{,} nonlinearities (cf.
for example, \cite{MR}, where the \anold{non}local CH equation with \anold{non}linear reaction term has been analyzed).
Furthermore, we note that interest in including reaction terms in CH-type systems is particularly growing due to the development of diffuse interface models of tumor growth that couple CH-type systems with nutrient diffusion and other equations. In fact, in this class of models\an{,} the parameter $\varphi$ represents the concentration of the tumor phase, and therefore, in this case, it is of particular interest not to have mass conservation, but to observe a possible growth or decrease of the tumor mass due to cell proliferation or death. In this regard, one can see, e.g., \anold{\cite{CRW, CGH, CGRS, CSS2, FGR, garcke, GARL_1, GLR, KS2, MRS,  RSchS,  SS, S}} and the
\an{references therein for results on well-posedness, long-\anold{time} behavior of solutions, asymptotic analyses and optimal control}.}

\betti{The solvency and existence of global and exponential attractors for the CHO equation with regular potential \eqref{regpot} were studied in \cite{M}, while in \cite{GGM} the authors investigated
     the case of the singular logarithmic potential \eqref{logpot}, establishing some regularization properties of the unique solution in finite time. Both in \an{two} and \an{three dimensions,} the existence of a global attractor has been proved\an{. Moreover, the authors of \cite{GGM}, exploiting the regularization effects they established in two dimensions and the associated strict separation property, proved the existence of an exponential attractor and the convergence to a single equilibrium.}}

In the present work, we couple \eqref{Isystem} with the following dynamic boundary condition for both $\mu$ and~$\phi$\an{. Namely, we} complement system \eqref{Isystem}~with
\anold{
\begin{align}
 \dt\phiG + \dn\mu - \DeltaG\muG = \gamma(\uG-\phiG),
	\quad
	\tau \dt\phiG + \dn\phi - \DeltaG\phiG + \FG'(\phiG) = \muG
 \quad  \hbox{on $\Sigma,$}
   \label{IdynBC}
\end{align}
with $\Sigma:=\an{\partial \Omega}\times (0,T)$ and $\Gamma:=\partial\Omega$, where} $\muG$ and $\phiG$ \anold{stand for} the traces of $\mu$ and~$\phi$, respectively,
$\DeltaG$~is the Laplace--Beltrami operator on the boundary,
$\uG$~is \anold{a boundary source/sink term that will later play the role of (boundary)} control \betti{on the boundary} $\Gamma:=\partial\Omega$ of the domain $\Omega$, $\FG'$ is the derivative of another \an{double-well shaped} potential~$\FG$, and $\dn$ denotes the outward normal derivative to $\Gamma$.
\betti{
Let us notice that here we consider the same viscosity coefficient $\tau$ in the bulk \anold{(cf. \eqref{Isystem})} and on the boundary \anold{(cf. \eqref{IdynBC})}  for simplicity.
\anold{H}owever, up to redefining the scalar product and the functional spaces, different values of $\tau$ could be encompassed in the analysis.}

\betti{\anold{T}his is the first time that the CHO equation is coupled with dynamic boundary conditions.
Indeed, in most works, system \eqref{Isystem} is endowed with \an{homogeneous} Neumann boundary conditions \anold{for both $\mu$ and $\phi$, that is,
\Beq
  \dn\mu = \dn\phi = 0 
  \quad \hbox{on $\Sigma$}.
  \non
\Eeq
\an{Those have, in the order, physical connection with the mass conservation, when no source terms are considered (i.e., when $\gamma=0$), and the contact angle problem as the no-flux condition for $\phi$ entails that the interface is orthogonal to the boundary.}}}

\betti{\anold{In the last decades}, physicists
have introduced the so-called dynamic boundary \anold{conditions}, in the sense that the kinetics, i.e., $\dt\phi$ appears explicitly in the boundary conditions, in order to account for the interaction of the
components with the walls for a confined system (see, e.g., \cite{FMD1, FMD2} and also \cite{Ketal, MMP} where numerical
simulations are performed).
Notice that for $\tau>0$, we obtain the viscous CH equation introduced in \cite{N}
(see \cite{BDS}, \cite{ES},  and \cite{MiZe} for the mathematical
analysis of the viscous CH equation with classical boundary conditions and \anold{\cite{GK}} for
dynamic boundary conditions and regular potentials). 
%
The interest in the analysis of CH-type equations coupled with dynamic boundary conditions 
has grown {more and more} in the recent literature. 
In this regard, we can quote {\cite{GiMiSchi} and} \cite{CGS2014}, 
where both the viscous and the \anold{non}viscous CH equations {with possibly singular potentials}, 
combined with these kinds of boundary conditions 
have been investigated
(see also \cite{GiMiSchi2} for the longtime behavior).
{In the former paper it is assumed that 
the bulk potential is dominant with respect to the boundary one,
while the assumption in the opposite direction is made in the latter, as well as in the present paper}.
Furthermore, we have to mention \anold{\cite{CF,CGS2017,CS,FN,GK,KS,KSnon,LW,MiZe}} and references therein, 
where problems related to the CH equation combined with different types of dynamic boundary conditions have been studied.}

\betti{Regarding the main scope of \anold{the manuscript}, since we are interested in the optimal control problem related to the previous PDE system coupled with dynamic boundary conditions, it is crucial to study the dependence of the solution on the \anold{variables}~$u$ and $\uG$ \anold{that will play the role of controls in the second part of the work}.
Indeed, a major part of the paper is devoted to the well-posedness of the system
and the continuous dependence of the solutions with respect to \anold{$u$} and $\uG$. 
These results can be proved under quite general assumptions on the potentials $F$ and $F_\Gamma$, requiring that the boundary potential $F_\Gamma$ is dominant with \anold{respect to} the bulk one $F$ (\anold{see}~also \cite{CGS2014} for similar assumptions). The case of the smooth  \eqref{regpot} and the singular \eqref{logpot} potentials can be included in the analysis, while \anold{the double obstacle potential} \eqref{obspot} \anold{can not be considered.}}

\betti{
Up to our knowledge, the only other reference dealing with optimal control for a CHO equation 
is \cite{CGRS4} where, however, \anold{the standard} Neumann homogeneous boundary conditions were considered for both $\phi$ and $\mu$.
Hence, well-posedness cannot be deduced from \cite{CGRS4} and the \anold{pre}existing literature.
Being our \anold{final aim} the \anold{application to} the \anold{aforementioned} control problem, we need first to prove well-posedness, regularity of solutions (including the separation property for $\phi$ in case of singular \anold{potentials}\an{)}, as well as the \frechet\ differentiability of the \anold{control to-state-mapping}. These are indeed preliminary results which will let us to prove \anold{first-order} necessary optimality conditions for the minimization problem  associated to the following \anold{tracking-type} cost functional
\an{\begin{align*}
   & \calJ(u,\uG;\phi,\phiG)
  := \frac {\alpha_1} 2 \intQ |\phi-\tarQ|^2
  + \frac {\alpha_2} 2 \intS |\phiG-\tarS|^2
  \non
  \\
 & \quad
  + \frac {\alpha_3} 2 \iO |\phi(T)-\anold{\tarO}|^2
  + \frac {\alpha_4} 2 \iG |\phiG(T)-\tarG|^2
  + \frac {\alpha_5} 2 \intQ |u|^2
  + \frac {\alpha_6} 2 \intS |\uG|^2,
 \end{align*}
where} $(u,\uG)$ belongs to a proper {set of admissible controls} (cf.~\eqref{defUad}) and $\phi$ and $\phiG$ are the components of the solution to problem \eqref{Isystem} and \eqref{IdynBC} (with suitable initial conditions) corresponding to the pair~$(u,\uG)$. }
\betti{\an{Above}, $\tarQ,$ $\tarS$, $\tarO$,  \anold{and}  $\tarG$ are given target functions, respectively \anold{in~$Q$ and $\Omega$, and on  $\Sigma$ and $\Gamma$, whereas} $\alpha_i$, $i=1,\dots,6$, are nonnegative constants.}

\betti{
\anold{After showing the existence of (at least) one optimal strategy, we investigate the first-order optimality conditions}.
\anold{In this direction, the main difficult} is to prove the \frechet\ differentiability 
of the control-to-state operator between suitable functional spaces,
which allows us to establish \anold{first-order} necessary optimality conditions 
in terms of the solution to a proper adjoint problem. The proof of solvability and uniqueness of solutions of the resulting adjoint system is one of the main difficulties of the paper, due to the fact that we are dealing with a CH-type equation with \anold{in}homogeneous mass source (CHO equation) as well as with dynamic boundary conditions \anold{of the same type}.}

\betti{
The paper is organized as follows. 
In the next section, we list our assumptions and \anold{notation}
and state our results.
The proofs of 
those regarding the well-posedness of the problem, 
the regularity and the continuous dependence of its solution on the \anold{variables $u$ and $\uG$} are given in Sections~\ref{EXIST-REG} and \ref{CONTDEP}.
Finally,  Section~\ref{CONTROL} is devoted to the analysis of
the control problem and the corresponding \anold{derivation of the} first-order necessary optimality conditions.}


\section{Statement of the problem and results}
\label{STATEMENT}
\setcounter{equation}{0}

In this section, we state precise assumptions and notations and present our results.
First of all, the set $\Omega\subset\erre^3$ 
is~assumed to be bounded, connected and smooth \an{(the lower
dimensional cases can be treated in a similar way)}.
As in the Introduction, $\dn$ and $\DeltaG$ stand for the outward normal derivative on $\Gamma:=\partial\Omega$
and the Laplace--Beltrami operator, respectively.
Furthermore, we denote by $\nablaG$ the surface gradient
and write $\misO$ and $\misG$ 
for the volume of $\Omega$ and the area of~$\Gamma$, respectively.
{Next, if $X$ is a Banach space, $\norma\cpto_X$ denotes its norm,
with the only exception for the norms in the $L^\infty$ spaces, 
which are denoted by~$\norma\cpto_\infty\,$. 
Moreover, the symbol $\norma\cpto_X$ will also stand for the norm in~$X^3$.
Similarly, when no confusion can arise, we simply write, e.g., $\L2X$ in place of $\L2{X^3}$.}
Moreover, for every Banach space~$X$, the symbols $X^*$ and $\<\cpto,\cpto>_X$ 
denote the dual space of $X$ and the dual pairing between $X^*$ and~$X$, respectively.
Furthermore, we introduce the \anold{shorthand}
\begin{align}
  & H := \Ldue \,, \quad  
  V := \Huno 
  \aand
  W := \Hdue,
  \label{defspaziO}
  \\
  & \HG := \LdueG \,, \quad 
  \VG := \HunoG 
  \aand
  \WG := \HdueG,
  \label{defspaziG}
  \\
  & \calH := H \times \HG \,, \quad
  \calV := \graffe{(v,\vG) \in V \times \VG : \ \vG = v\suG}
  \non
  \\
  & \aand
  \calW := \bigl( W \times \WG \bigr) \cap \calV \,.
  \label{defspaziprod}
\end{align}
Finally, we define
\Beq
  \mean(z,\zG) = \frac { \iO z + \iG \zG } { \misO + \misG }
  \quad \hbox{for $(z,\zG)\in\calH$}
  \label{defmean}
\Eeq
and term it \anold{extended} {mean value} of~$(z,\zG)$.

{To help the reader, we adopt the subscript $\Gamma$ in pairs of type $(v,\vG)$
only when the second component is the trace of the first one,
while we use superscripts in the opposite case (as~in \eqref{defmean})}.

\vskip 2mm

Now, we list our assumptions.
As for the structure of our system, 
we postulate~that
\begin{align}
  & \hbox{$\tau$ and $\gamma$ are positive real numbers}
  \label{hptg}
  \\
  & F,\FG : \erre \to (-\infty,+\infty]
  \quad \hbox{admit the decompositions}
  \non
  \\
  & \quad F = \Beta + \Pi
  \aand
  \FG = \BetaG + \PiG
  \quad \hbox{where}
  \label{hppot}
  \\
  & \Beta,\, \BetaG : \erre \to [0,+\infty]
  \quad \hbox{are convex and l.s.c.\ with} \quad
  \Beta(0) = \BetaG(0) = 0
  \qquad
  \label{hpBeta}
  \\
  \separa
  & \Pi,\, \PiG : \erre \to \erre
  \quad \hbox{are of class $C^2$ with \Lip\ continuous derivatives}.
  \qquad
  \label{hpPi}
\end{align}
We set, for convenience,
\Beq
  \beta := \partial\Beta \,, \quad
  \betaG := \partial\BetaG \,, \quad
  \pi := \Pi'
  \aand
  \piG := \PiG'
  \label{defbetapi}  
\Eeq
and require the compatibility condition
\begin{align}
  & D(\betaG) \subseteq D(\beta)
  \aand
  |\betaz(r)| \leq C^* \bigl( |\betaGz(r)| + 1 \bigr)
  \non
  \\
  & \quad \hbox{for some $C^*>0$ and every $r\in D(\betaG)$} \,.
  \label{hpCC}
\end{align}
\Accorpa\HPstruttura hptg hpCC
\rev{In the above lines, the symbol $\partial$ indicates the subdifferential operator: we refer to \cite{Brezis} for more details.}
{Here}, the symbols $D(\beta)$ and $D(\betaG)$ 
denote the domains of $\beta$ and~$\betaG$, respectively.
More generally, we use the notation $D(\sigma)$ 
for every maximal monotone graph $\sigma$ in $\erre\times\erre$.
Moreover, for $r\in D(\sigma)$,
$\sigma^\circ(r)$ stands for the element of $\sigma(r)$ having minimum modulus.
Finally, we still term $\sigma$ the maximal monotone operators induced on $L^2$ spaces.

{%
\Brem
Condition \eqref{hpCC} prescribes that the boundary potential is dominant compared to the one in the bulk.
This condition has been assumed in several papers (see \cite{CGRS4} among others).
A~condition in the opposite direction has instead been postulated in \cite{GiMiSchi}.
\Erem
}%

For the data, we make the assumptions listed below.
They involve the constant $M$ that is related to the optimal control problem we state later~on.
At the present \anold{stage}, one can think that $M$ is just a given constant.
We assume~that
\begin{align}
  & u \in \LQ\infty, \quad
  \uG \in \LS\infty
  \quad \hbox{with} \quad
  \norma u_\infty \leq M
  \aand
  \norma\uG_\infty \leq M
  \label{hpu}
  \\
  & (\phiz \,,\, \phiGz) \in \calV \,, \quad
  \Beta(\phiz) \in \Luno
  \aand
  \BetaG(\phiGz) \in \LunoG 
  \qquad
  \label{hpphiz}
  \\
  & \hbox{if} \quad \mz := \mean(\phiz,\phiGz)
  \aand 
  \rho := \frac{ M} \gamma \,
  \quad \hbox{then} 
  \non
  \\
  & \quad \hbox{$-\mz^- - \rho$ \ and \ $\mz^+ + \rho$ \ belong to \ $\mathop{\rm int}\nolimits D(\betaG)$}
  \label{hpmz}
\end{align}
\Accorpa\HPdati hpu hpmz
\gianni{where $\mz^-$ and $\mz^+$ denote the negative and positive part of~$\mz$, respectively.}

Let us {now} come to our notion of solution.
It corresponds to the weak form of the problem presented in the Introduction
that one deduces this way:
given an arbitrary $(v,\vG)\in\calV$,
one formally multiplies the equations {in} \eqref{Isystem} by $v$
and those {in} \eqref{IdynBC} by $\vG$ and integrates by parts in space.
Then, the normal derivatives $\dn\mu$ and $\dn\phi$ {simplify}.
However, although a lower level of regularity would be sufficient to make the definition meaningful,
we require that the solution is smoother.
Thus, we look for a {six}-tuple $\soluz$ (or~a triplet of pairs) 
that satisfies the regularity requirements
\begin{align}
  & (\mu,\muG) \in \L2\calW 
  \label{regmu}
  \\
  & (\phi,\phiG) \in \H1\calH \cap \L\infty\calV \cap \L2\calW
  \label{regphi}
  \\
  & (\xi,\xiG) \in \L2\calH 
  \label{regxi}
\end{align}
\Accorpa\Regsoluz regmu regxi
and solves
\begin{align}
  & \iO \dt\phi \, v
  + \iG \dt\phiG \, \vG
  + \iO \nabla\mu \cdot \nabla v
  + \iG \nablaG\muG \cdot \nablaG\vG
  \non
  \\
  & = \gamma \iO (u-\phi) v
  + \gamma \iG (\uG-\phiG) \vG
  \non
  \\
  & \quad \hbox{\aet\ and for every $(v,\vG)\in\calV$},
  \label{prima}
  \\
  \separa
  & \tau \iO \dt\phi \, v
  + \tau \iG \dt\phiG \, \vG
  + \iO \nabla\phi \cdot \nabla v
  + \iG \nablaG\phiG \cdot \nablaG\vG
  \non
  \\
  & \quad {}
  + \iO \bigl( \xi + \pi(\phi) \bigr) v
  + \iG \bigl( \xiG + \piG(\phiG) \bigr) \vG
  = \iO \mu v 
  + \iG \muG \vG
  \non
  \\
  & \quad \hbox{\aet\ and for every $(v,\vG)\in\calV$},
  \label{seconda}
  \\
  & \xi \in \beta(\phi) \quad \aeQ
  \aand
  \xiG \in \betaG(\phiG) \quad \aeS
  \label{terza}
  \\
  & \phi(0) = \phiz
  \quad \aeO \,.
  \label{cauchy}
\end{align}
\Accorpa\Pbl prima cauchy
{We did not \an{explicitly} assume that $\phiG(0)=\phiGz$ \aeG\
since this condition is clearly implied by~\eqref{cauchy}.}

\Brem
\label{Strong}
One can show that every solution to \Pbl\ in the sense specified above
actually solves the problem in a strong form, i.e., 
equations \eqref{Isystem} and \eqref{IdynBC} are satisfied
with $F'(\phi)$ and $\FG'(\phiG)$ replaced by $\xi+\pi(\phi)$ and $\xiG+\piG(\phiG)$, respectively.
This could be proved by suitably applying the forthcoming Lemma~\ref{Elliptic}.
\Erem

The first results of ours regard the existence of a solution
and a stability estimate.

\Bthm
\label{Existence}
Assume \HPstruttura\ on the structure of the system and \HPdati\ on the data.
Then, problem \Pbl\ has at least one solution $\soluz$ satisfying the regularity requirements \Regsoluz\
as well as the stability estimate
\Beq
  {\norma{(\mu,\muG)}_{\L2\calW}
  + \norma{(\phi,\phiG)}_{\H1\calH\cap\L\infty\calV\cap\L2\calW}
  + \norma{(\xi,\xiG)}_{\L2\calH}}
  \leq K_1
  \label{stability}
\Eeq
with a {positive} constant $K_1$ that depends only on the structure of the system,
$\Omega$, $T$, the constant $M$ {appearing} in \eqref{hpu} 
and the initial datum~$\phiz$.
In particular, $K_1$~{is independent of} the pair $(u,\uG)$.
\Ethm

Under further regularity assumptions {on the data}, we can prove the existence of a more regular solution.

\Bthm
\label{Regularity}
Let assumptions \HPstruttura\ on the structure hold and, in addition to \HPdati\ for the data,
suppose that $(u,\uG)$ and $\phiz$ satisfy for some \anold{positive} constant~$M'$
\begin{align}
  & (u,\uG) \in \H1\calH
  \quad \hbox{with} \quad
  {\max\graffe{\norma{\dt u}_{\L2H},\norma{\dt\uG}_{\L2\HG}} \leq M'}
  \label{hpureg}
  \\
  & \anold{(\phiz,\phiGz)} \in \calW 
  \aand
  (\betaz(\phiz),\betaGz(\phiGz)) \in \calH .
  \label{hpphizreg}
\end{align}
\Accorpa\HPdatireg hpureg hpphizreg
Then, problem \Pbl\ has at least a solution $\soluz$ 
that also satisfies
\begin{align}
  & (\mu,\muG) \in \L\infty\calW \,, \quad
  (\phi,\phiG) \in \W{1,\infty}\calH \cap \H1\calV \cap \L\infty\calW
  \qquad
  \non
  \\
  & \aand (\xi,\xiG) \in \L\infty\calH 
  \label{regularity}  
  \\
  & \norma{(\mu,\muG)}_{\anold{\L\infty\calW}}
  + \norma{(\phi,\phiG)}_{\W{1,\infty}\calH\cap\H1\calV\cap\L\infty\calW}
  \non
  \\
  & \quad {}
  + \norma{(\xi,\xiG)}_{\L\infty\calH}
  \leq K_2
  \qquad
  \label{stabilitybis}
\end{align}
with a \anold{positive} constant $K_2$ that depends only on the structure of the system, 
$\Omega$, $T$, the constants $M$ and $M'$ and initial datum~$\phiz$.
{In particular, as a consequence of the continuous embedding $\calW\emb\Linfty\times\LxG\infty$,}
the first four components of the solution are {uniformly} bounded
and their $L^\infty$ norms are estimated by a constant {with the same dependence of}~$K_2$.
\Ethm

In the case of everywhere defined potentials {(e.g., \eqref{regpot})}, the boundedness of the component $\phi$
already given by Theorem~\ref{Regularity} is satisfactory.
If instead potentials of logarithmic type are considered,
one essentially deals with smooth potentials only if the values of the component $\phi$
are far away from the end-points of the domain.
Our last results regarding the properties of the solutions \an{cover} this {scenario}.
However, {regular potentials are also included in our analysis}.
Indeed, we \anold{require}~that
\begin{align}
  & \hbox{either \ $D:=D(\beta)=D(\betaG)=\erre$ \ or \ $D:=D(\beta)=D(\betaG)=(-1,1)$}
  \label{hpD}
  \\
  & \hbox{$\beta,\,\betaG:D\to\erre$ \ are single-valued and locally \Lip\ continuous}.
  \label{hpF}
\end{align}
Due to the last assumption, we can {ignore} the components $\xi$ and $\xiG$ {from the six-tuple} $\soluz$
and consider just the quadruplet $\sol$ when speaking of a solution,
{as $\xi=\beta(\phi)$ and $\xiG=\betaG(\phiG)$}.

\Bthm
\label{Separation}
In addition to assumptions \HPstruttura\ on the structure,
assume that \accorpa{hpD}{hpF} are satisfied.
Moreover, assume that the data satisfy \HPdati\ and \HPdatireg\ as well~as
\Beq
  \phiz\in D_0 \quad \aeO
  \quad \hbox{for some compact subset $D_0\subset D$}.  
  \label{hpsepar}
\Eeq
Then, \anold{there} exists a solution $\sol$ to problem \Pbl\ satisfying 
\Beq
  \phimin \leq \phi \leq \phimax
  \quad \aeQ
  {\aand
  \phimin \leq \phiG \leq \phimax
  \quad \aeS}
  \label{separation}
\Eeq
for some constants $\phimin,\phimax\in D$
that depend only on the structure of the system, 
$\Omega$, $T$, the constants $M$ and $M'$ and the initial datum \an{$\phiz$}.
\Ethm

The next result regards 
{the uniqueness of the solution and its continuous dependence on the pair $(u,\uG)$}.
We are not able to prove uniqueness in the general case, unfortunately.
However, the result we state below covers both the case of classical and logarithmic potentials.

\Bthm
\label{Contdep}
Under the same assumptions of Theorem~\ref{Separation} on the structure and the data, 
the solution to problem \Pbl\ satisfying the theses of Theorems~\ref{Regularity} and~\ref{Separation} is unique.
Moreover, if $(u_i,u_{\Gamma,i})$, $i=1,2$, are two choices of $(u,\uG)$
{and $(\mu_i,\mu_{\Gamma,i},\phi_i,\phi_{\Gamma,i})$ are the corresponding solutions,}
the following estimate
{%
\begin{align}
  & \norma{(\mu_1,\mu_{\Gamma,1})-(\mu_2,\mu_{\Gamma,2})}_{\L2\calV}
  + \norma{(\phi_1,\phi_{\Gamma,1})-(\phi_2,\phi_{\Gamma,2})}_{\H1\calH\cap\L\infty\calV}
  \non
  \\
  & \leq K_3 \, \norma{(u_1,u_{\Gamma,1})-(u_2,u_{\Gamma,2})}_{\L2\calH}
  \label{contdep}
\end{align}
}%
holds true with a constant $K_3$ that depends only on the structure of the system,
$\Omega$, $T$, the constants $M$ and $M'$ and the initial datum~$\phiz$.
\Ethm

Once the well-posedness of the problem and some regularity of its solution are established,
one can deal with {a corresponding control problem}.
Even though a {slightly} more general cost functional could be considered
(see the forthcoming Remark~\ref{Moregeneral}),
the problem we aim at studying is the following:
\begin{align}
  & \hbox{\sl Minimize the cost functional}
  \non
  \\
  & \calJ(u,\uG;\phi,\phiG)
  := \frac {\alpha_1} 2 \intQ |\phi-\tarQ|^2
  + \frac {\alpha_2} 2 \intS |\phiG-\tarS|^2
  \non
  \\
  & \phantom{\calJ(u,\uG;\phi,\phiG):=}
  + \frac {\alpha_3} 2 \iO |\phi(T)-\anold{\tarO}|^2
  + \frac {\alpha_4} 2 \iG |\phiG(T)-\tarG|^2
  \non
  \\
  & \phantom{\calJ(u,\uG;\phi,\phiG):=}
  + \frac {\alpha_5} 2 \intQ |u|^2
  + \frac {\alpha_6} 2 \intS |\uG|^2
  \label{cost}
  \\
  \separa
  & \hbox{\sl {subject to} the constraints}
  \non
  \\
  & i)\quad \hbox{\sl $(u,\uG)$ belongs to the {set of admissible controls}}
  \non
  \\
  & \hskip 4em \Uad := \bigl\{(u,\uG)\in{\bigl(\LQ\infty\times\LS\infty\bigr)\cap\H1\calH}:
  \non
  \\
  & \hskip 8em \umin\leq u\leq\umax \ \ \hbox{and} \ \ \uGmin\leq\uG\leq\uGmax \,,
  \non
  \\
  & \hskip 8em \norma{\dt u}_{\L2H}\leq M' \ \ \hbox{and} \ \ \norma{\dt\uG}_{\L2\HG}\leq M'\bigr\}
  \label{defUad}
  \\
  & ii) \quad \hbox{\sl $\phi$ and $\phiG$ are the components of the solution $\sol$}
  \non
  \\
  & \qquad \hbox{\sl  to problem \Pbl\ corresponding to the pair~$(u,\uG)$.}
  \label{cts}
\end{align}
\Accorpa\ControlPbl cost cts
For the quantities appearing in \eqref{cost} and \eqref{defUad}\anold{,} we require~that
\begin{align}
  & \alpha_i \enskip (i=1,\dots,6) \enskip \hbox{are nonnegative constants}
  \label{coeff}
  \\
  & \tarQ\in\LQ2 , \quad
  \tarS\in\LS2 , \quad
  \tarO\in\Ldue
  \aand
  \tarG\in\LdueG
  \label{target} 
  \\
  & \umin,\,\umax \in \LQ\infty
  \aand
  \uGmin,\,\uGmax\in\LS\infty
  \label{hpUad}
  \\
  & \hbox{condition \eqref{hpmz} is satisfied with \ }
  \non
  \\
  & \quad
  M := \max\graffe{\norma\umin_\infty\,,\,\norma\umax_\infty\,,\,\norma\uGmin_\infty\,,\,\norma\uGmax_\infty}
  \label{compatUad}
  \\
  & \hbox{$M'$ is a positive constant}.
  \label{emmeprimo}
\end{align}
Finally, we assume that
\Beq
  \hbox{$\Uad$ is nonempty}.
  \label{nonempty}
\Eeq
\Accorpa\HPcontrol coeff nonempty

{%
\Brem
\label{Remcontrol}
Notice that \eqref{nonempty} implies that 
$\umin\leq\umax$ \aeQ\ and $\uGmin\leq\uGmax$ \aeS.
On the contrary, these conditions do not imply \eqref{nonempty},
because of the restriction on the time derivative.
We also observe that the special case $\umin=\umax$ is allowed
provided that \eqref{nonempty} still holds.
In this case, the only variable from the pair $(u,\uG)$ is $\uG$
and problem \ControlPbl\ becomes a boundary control problem.
Accordingly, the bulk integral in the necessary condition of optimality we describe below disappears
and just the boundary integral survives.
Similarly, we can consider \an{just} the \an{distributed} control problem with the control acting only in the bulk.
\Erem
}%

Thanks to the well-posedness and regularity results stated above,
we can assume that problem \Pbl\ is well-posed and that its solution is smooth
and satisfies the proper uniform bounds as well as the separation property.
As for the latter, we recall that 
the main assumptions we have made on the structure were~\accorpa{hpD}{hpF}.
However, to develop the complete {control} theory, we need more regularity on the potentials.
Therefore, we reinforce \eqref{hpF} by requiring~that
\Beq
  \hbox{$F,\,\FG:D\to\erre$ \ are functions of class $C^3$}.
  \label{hpFreg}
\Eeq
{Notice that this condition is still met by the regular and logarithmic potentials \eqref{regpot} and~\eqref{logpot}}.
Here is our first result,
in which it is understood that all the assumption\an{s} we have listed are satisfied.

\Bthm
\label{OKcontrol}
The control problem \ControlPbl\ has at least one solution $(\ustar,\uGstar)$.
\Ethm

\betti{Our last result consists in finding} necessary conditions for optimality.
We will prove the following: 
if $(\ustar,\uGstar)$ is an optimal control and $\solstar$ is the corresponding state, 
then, we have that
\Beq
  {\intQ (\gamma p + \alpha_5 \ustar) (u-\ustar)
  + \intS (\gamma\pG + \alpha_6 \uGstar) (\uG-\uGstar)}
  \geq 0
  \quad \hbox{for every $(u,\uG)\in\Uad$}
  \non
\Eeq
where $p$ and $\pG$ are the components of the solution $\soluza$ to a proper adjoint problem
which is discussed in detail in the last section.
Here, we {only anticipate} that a strong form of it is the following backward {in time} problem,
$\pG$~and $\qG$ being the traces of $p$ and $q$ on the boundary:
\Bsist
  && - \dt(p+\tau q) - \Delta q + F''(\phistar) \, q {{}+\gamma p}
  = \alpha_1 (\phistar-\tarQ)
  \non
  \\
  && \aand
  -\Delta p = q
  \quad \hbox{in $Q$}
  \non
  \\
  && -\dt(\pG+\tau\qG) + \dn q - \DeltaG\qG + \FG''(\phiGstar) \, \qG {{}+\gamma\pG}
  = \alpha_2 (\phiGstar-\tarS)
  \non
  \\
  && \aand
  \dn p - \DeltaG\pG = \qG
  \quad \hbox{on $\Sigma$}
  \qquad  
  \non
  \\
  && (p+\tau q)(T) = \alpha_3(\phistar(T)-\tarO)
  \quad \hbox{in $\Omega$}
  \non
  \\
  && \aand
  (\pG+\tau\qG)(T) = \alpha_4(\phiGstar(T)-\tarG)
  \quad \hbox{on $\Gamma$} 
  \non
\Esist
{where $\phistar$ and $\phiGstar$ are the components of the solution $\solstar$
corresponding to the optimal control $(\ustar,\uGstar)$}.
However, this is completely formal without \anold{specific} assumptions on the ingredients of the cost functional.

The rest of the paper is organized as follows.
We continue the present section by fixing some notations and recalling some tools.
The next Section~\ref{EXIST-REG} is devoted to the proof of 
our existence and regularity results,
{whereas} continuous dependence and uniqueness {are} proved in Section~\ref{CONTDEP}.
Finally, the control problem is discussed in the last section.

\medskip

Throughout the paper, we will repeatedly use the Young inequality
\Beq
  a\,b \leq \delta\,a^2 + \frac 1{4\delta} \, b^2
  \quad \hbox{for all $a,b\in\erre$ and $\delta>0$}
  \label{young}
\Eeq
as well as the \Holder\ and Schwarz inequalities.
Moreover, we employ the abbreviations
\Beq
  Q_t := \Omega\times(0,t)
  \aand
  \Sigma_t := \Gamma\times(0,t)
  \quad \hbox{for $t\in(0,T]$}.
  \label{defQtSt}
\Eeq
Next, we introduce some tools related to the notion of generalized mean value given by~\eqref{defmean}.
We observe that
\Beq
  \iO (v-\mean(v,\vG)) 
  + \iG (\vG-\mean(v,\vG))
  = 0
  \quad \hbox{for every $(v,\vG)\in\calV$}.
  \label{meanzero}
\Eeq
It \anold{then} follows that
\Beq
  \calV 
  = \calVz \oplus \Span\graffe{(1,1)}
  \label{decomp}
\Eeq
where 
\Beq
  \calVz := \graffe{ (v,\vG)\in\calV : \ \mean(v,\vG) = 0 }.
  \label{defVz}
\Eeq
Of course, the components of the pair $(1,1)$ in \eqref{decomp}
are the constant functions~$1$ on $\Omega$ and~$\Gamma$, respectively.
More generally, if $y$ is any real number {(e.g., some mean value)}, the same symbol $y$
also denotes the corresponding constant functions on $\Omega$ and~$\Gamma$.
The same convention is used if $y$ is time dependent.
One also can deduce from \eqref{decomp} the Poincar\'e type inequality
\Beq
  \normaHH{(v,\vG)}^2
  \leq \CO \bigl( \normaH{\nabla v}^2 + \normaHG{\nablaG\vG}^2 + |\mean(v,\vG)|^2 \bigr)
  \quad \hbox{for every $(v,\vG)\in\calV$}
  \label{poincare}
\Eeq
where $\CO$ depends only on~$\Omega$.
It follows that the {map}
\Beq
  \calV \ni (v,\vG) \an{\mapsto}
  \bigl( \normaH{\nabla v}^2
  + \normaHG{\nablaG\vG}^2
  + |\mean(v,\vG)|^2 \bigr)^{1/2}
  \non
\Eeq
yields a Hilbert norm on $\calV$ that is equivalent to the natural one.
Finally, we will take advantage of the following \anold{regularity} result
{(that can be found in~\cite[Lem.~3.1]{CGS13}).}

\Blem
\label{Elliptic}
Let $\sigma:\erre\to\erre$ be monotone and \Lip\ continuous {function}
and assume that $(w,\wG)\in\calV$ and $(g,\gG)\in\calH$ verify
\Beq
  \iO \nabla w \cdot \nabla v
  + \iG \nabla\wG \cdot \nabla\vG
  + \iG \sigma(\wG) \vG
  = \iO g v
  + \iG \gG \vG
  \quad \hbox{for every $(v,\vG)\in\calV$}.
  \label{elliptic}
\Eeq
Then, we have
\Beq
  (w,\wG) \in \calW
  \aand
  \normaWW{(w,\wG)}
  + \normaHG{\sigma(\wG)}
  \leq \CO \bigl( \normaVV{(w,\wG)} + \normaHH{(g,\gG)} \bigr)
  \label{ellreg}
\Eeq
where $\CO$ depends only on~$\Omega$ and not on~$\sigma$.
Moreover, $(w,\wG)$ solves the boundary value problem
\Beq
  - \Delta w = g
  \quad \aeO
  \aand
  \dn w - \DeltaG\wG + \an{\sigma}(\wG) = \gG
  \quad \aeG \,.
  \label{bvpbl}
\Eeq
\Elem

We conclude this section by stating a general rule 
concerning the constants that appear in the estimates we perform in the following.
We use the small-case symbol $\,c\,$ for a generic constant
whose actual values may change from line to line and even within the same line
and depend only on~$\Omega$, the structure of the system,
and the constants and the norms of the functions involved in the assumptions of the statements.
In particular, the \anold{value} of $\,c\,$ may depend on the constant $M$ that appears in \eqref{hpu}
but they are independent of $u$ and~$\uG$
and it is also clear that they do not depend on the parameter $\,\eps\,$ we introduce in Section~\ref{EXISTENCE}.
A~small-case symbol with a subscript like $c_\delta$
indicates that the constant may depend on the parameter~$\delta$, in addition.
On the contrary, we mark precise constants that we can refer~to
by using different symbols 
(e.g.,~capital or \gianni{Greek} letters like in \eqref{hpCC} and~\eqref{hpmz}).


\section{Existence and regularity}
\label{EXIST-REG}
\setcounter{equation}{0}

In this section, we prove Theorems~\ref{Existence}, \ref{Regularity} and~\ref{Separation} \anold{whose proofs} are inspired by the papers \cite{CGRS4} and~\cite{CGS13}.
\an{In this direction}, \an{as a first step,} we introduce \anold{a} regularized problem
that is obtained by replacing the functionals $\Beta$ and $\BetaG$ and the operators $\beta$ and $\betaG$
by their Moreau--Yosida regularizations 
$\Betaeps$, $\BetaGeps$, $\betaeps$ and~$\betaGeps$, respectively, {with $\eps\in(0,1)$}
(see, e.g., \cite[pp.~28 and~39]{Brezis}).
Thus, \gianni{the regularized problem} consists in finding a quadruplet $\soluzeps$ 
(or~a pair of pairs) 
that satisfies the regularity requirements
\begin{align}
  & (\mueps,\muGeps) \in \L2\calW 
  \label{regmueps}
  \\
  & (\phieps,\phiGeps) \in \H1\calH \cap \L\infty\calV \cap \L2\calW
  \label{regphieps}
\end{align}
\Accorpa\Regsoluzeps regmueps regphieps
and solves
\begin{align}
  & \iO \dt\phieps \, v
  + \iG \dt\phiGeps \, \vG
  + \iO \nabla\mueps \cdot \nabla v
  + \iG \nablaG\muGeps \cdot \nablaG\vG
  \non
  \\
  & = \gamma \iO (u-\phieps) v
  + \gamma \iG (\uG-\phiGeps) \vG
  \non
  \\
  & \quad \hbox{\aet\ and for every $(v,\vG)\in\calV$}
  \label{primaeps}
  \\
  \separa
  & \tau \iO \dt\phieps \, v
  + \tau \iG \dt\phiGeps \, \vG
  + \iO \nabla\phieps \cdot \nabla v
  + \iG \nablaG\phiGeps \cdot \nablaG\vG
  \non
  \\
  & \quad {}
  + \iO (\betaeps + \pi)(\phieps) \, v
  + \iG (\betaGeps + \piG)(\phiGeps) \, \vG
  = \iO \mueps v 
  + \iG \muGeps \vG
  \non
  \\
  & \quad \hbox{\aet\ and for every $(v,\vG)\in\calV$},
  \label{secondaeps}
  \\
  & \phieps(0) = \phiz
  \quad \aeO \,.
  \label{cauchyeps}
\end{align}
\Accorpa\Pbleps primaeps cauchyeps
\anold{Let us claim that this problem} has a unique solution \anold{for every $\eps \in (0,1)$}.
We do not prove this fact and just make two \anold{remarks}.
First, uniqueness is not used in the sequel.
Nevertheless, it could be proved as we do for Theorem~\ref{Contdep}
since $\an{\betaeps}$ and $\betaGeps$ are \Lip\ continuous.
As for existence, its proof could be obtained by discretizing the problem
by means of a Faedo--Galerkin scheme as done\an{, e.g.,} in~\cite{CGS13}.
The estimates we formally establish in the sequel 
can be performed in a rigorous way on the discrete solution
since, e.g., even its time derivative actually is an admissible test function.
So, we start from the solution $\soluzeps$ to \Pbleps\ and perform some formal estimates
that allow us to let $\eps$ tend to zero.
To this end, we make some preliminary observations.

We notice that \eqref{hpBeta} and \eqref{defbetapi} hold for the approximation.
More precisely, \anold{for every $\eps \in (0,1)$,} we have~that
\begin{align}
  & 0 \leq \Betaeps(r) \leq \Beta(r)
  \aand
  0 \leq \BetaGeps(r) \leq \BetaG(r)
  \quad \hbox{for every $r\in\erre$}
  \label{propBetaeps}
  \\
  & |\betaeps(r)| \leq |\betaz(r)|
  \aand
  |\betaGeps(r)| \leq |\betaGz(r)|
  \quad \hbox{for every $r\in D(\beta)$} .
  \label{propbetaeps}
\end{align}
Furthermore, \eqref{hpCC}~also holds true for $\betaeps$ and $\betaGeps$
with a similar constant that does not depend on~$\eps$ (see~\cite[Lemma 4.4]{CaCo}).
We thus write
\Beq
  |\betaeps(r)| \leq C^* \bigl( |\betaGeps(r)| + 1 \bigr)
  \quad \hbox{for every $r\in\erre$}
  \label{propCCeps}
\Eeq
with the same constant $C^*$, without loss of generality.
We also deduce~that
\Beq
  |\Betaeps(r)| \leq C^* \bigl( |\BetaGeps(r)| + |r| \bigr)
  \quad \hbox{for every $r\in\erre$} 
  \label{dapropCCeps}
\Eeq
just by integration.
Next, since $\betaeps$ and $\betaGeps$ have the same sign,
we see that \eqref{propCCeps} and the Young inequality imply
\Beq
  \betaGeps(r) \betaeps(r)
  \geq \frac 1 {2C^*} \, |\betaeps(r)|^2 - C_*
  \quad \hbox{for every $r\in\erre$}
  \label{prodbetaeps}
\Eeq
with a similar constant~$C_*$.
We also notice that the \anold{domains} inclusion $D(\betaG)\subseteq D(\beta)$ 
\anold{prescribed in~\eqref{hpCC}} and \eqref{hpmz} imply the existence of positive constants $\delta_0$ and $C_0$ such~that
\begin{align}
  & \betaeps(r) (r-r_0)
  \geq \delta_0 |\betaeps(r)| - C_0
  \aand
  \betaGeps(r) (r-r_0)
  \geq \delta_0 |\betaGeps(r)| - C_0
  \non
  \\
  & \quad \hbox{for every $r\in\erre$, $r_0\in[-\mz^- -\rho,\mz^+ +\rho]$ and $\eps\in(0,1)$}.
  \label{trickMZ}
\end{align}
This is a generalization of \cite[Appendix, Prop. A.1]{MiZe}.
The detailed proof given in \cite[p.~908]{GiMiSchi} with a fixed $r_0$ 
also works in the present case with minor changes.


\subsection{Existence}
\label{EXISTENCE}

We {now} prove Theorem~\ref{Existence}.
We start from the solution $\soluzeps$ to the approximating problem \Pbleps\
and perform a number of a priori estimates {that are uniform with respect to~$\eps$}.
The first bound we establish regards the {mean value} of~$(\phieps,\phiGeps)$.

\step
Estimate of a mean value

{We set} for brevity
\Beq  
  m := \mean(\phieps,\phiGeps)
  \aand
  \omega := \mean(u,\uG) 
  \quad \aet
  \label{defm}
\Eeq
\betti{and} \anold{test  \eqref{primaeps}  evaluated  at the time $t$ by $(1,1)/(\misO+\misG)$ to obtain that}
\Beq
  m'(t) + \gamma \, m(t) = \gamma \, \omega(t)
  \quad \aet.
\Eeq
On the other hand, \eqref{cauchyeps} implies that
$m(0)=\mz:=\mean(\phiz,\phiGz)$.
Therefore, $m$ does not depend on $\eps$ and we have~that
\Beq
  m(t) = \mz \, e^{-\gamma t} + \iot e^{-\gamma (t-s)} \omega(s) \, ds
  \quad \hbox{for every $t\in[0,T]$}.
  \non
\Eeq
\gianni{Now, for every $t\in[0,T]$, it results that
\begin{align}
  & -\mz^-
  \leq -\mz^- \, e^{-\gamma t}  
  \leq \mz \, e^{-\gamma t}
  \leq \mz^+ \, e^{-\gamma t}
  \leq \mz^+
  \non
  \\
  & \Bigl| \iot e^{-\gamma (t-s)} \omega(s) \, ds \Bigr|
  \leq M \Bigl| \iot e^{-\gamma (t-s)} \, ds \Bigr|
  \leq \frac M\gamma 
  = \rho 
  \non
\end{align}
the latter since $|\omega|\leq M$ \aet\ by \eqref{hpu}
and the definition of~$\rho$ (see \eqref{hpmz}).
Hence, we conclude~that
\Beq
  \gianni{-\mz^-} -\rho \leq \mean(\phieps(t),\phiGeps(t)) \leq \gianni{\mz^+} + \rho
  \quad \hbox{for every $t\in[0,T]$ and $\eps\in(0,1)$}.
  \label{stimamean}
\Eeq
}%

\step
An auxiliary estimate

We set for convenience
\Beq
  \aeps := \mean(\mueps,\muGeps) 
  \label{defalpha}
\Eeq
and still use the notation $m$ introduced in \eqref{defm}.
We \an{then} test \eqref{secondaeps} by $(\phieps-m,\phiGeps-m)$ and obtain \aet
\Bsist
  && \iO \betaeps(\phieps) (\phieps-m)
  + \iG \betaGeps\an{(\phiGeps)}(\phiGeps-m)
  + \iO |\nabla\phieps|^2
  + \iG |\nablaG\phiGeps|^2
  \non
  \\
  && = - \tau \iO \dt\phieps (\phieps-m)
  - \tau \iG \dt\phiGeps (\phiGeps-m)
  \non
  \\
  && \quad {}
  - \iO \pi(\phieps) (\phieps-m)
  - \iG \piG(\phiGeps) (\phiGeps-m)
  \non
  \\
  && \quad {}
  + \iO (\mueps-\aeps)(\phieps-m)  
  + \iG (\muGeps-\aeps)(\phiGeps-m) 
  \non
\Esist
{where, in the last line, we have subtracted 
$\iO\aeps(\phieps-m)+\iG\aeps(\phiGeps-m)$, which is zero due to~\eqref{meanzero}}.
Now, we owe to \eqref{stimamean} and apply~\eqref{trickMZ}.
\anold{Schwarz's} inequality in~$\calH$ 
and \anold{the} \Lip\ continuity of $\pi$ and~$\piG$  \anold{imply} that
\begin{align}
  & \delta_0 \iO |\betaeps(\phieps)|
  + \delta_0 \iG |\betaGeps(\phiGeps)|
  - 2 \, C_0
  \non
  \\
  & \leq \tau \normaHH{(\dt\phieps,\dt\phiGeps)} \, \normaHH{(\phieps-m,\phiGeps-m)}
  \non
  \\
  & \quad {}
  + c \, \bigl( \normaHH{(\phieps,\phiGeps)} + 1 \bigr) \normaHH{(\phieps-m,\phiGeps-m)}
  \non
  \\
  & \quad {}
  + \normaHH{(\mueps-\aeps,\muGeps-\aeps)} \, \normaHH{(\phieps-m,\phiGeps-m)} \,.
  \non
\end{align}
On the other hand, by testing \eqref{secondaeps} by $(1,1)$,
we obtain \aet\ that
\Bsist
  && (\misO+\misG) \, \aeps
  = \iO \betaeps(\phieps)
  + \iG \betaGeps(\phiGeps)
  \non
  \\
  && \quad {}
  + \tau \iO \dt\phieps
  + \tau \iG \dt\phiGeps
  + \iO \pi(\phieps)
  + \iG \piG(\phiGeps) .
  \non
\Esist
By combining {the above equalities and inequalities and owing \an{once again} to} \eqref{stimamean}, we conclude~that 
\begin{align}
  & |\aeps|
  \leq c \, \bigl\{
    \normaHH{(\dt\phieps,\dt\phiGeps)} \, \bigl( \normaHH{(\phieps,\phiGeps)} + 1 \bigr)
  \non
  \\
  & \quad {}
  + \normaHH{(\mueps-\aeps,\muGeps-\aeps)} \, \bigl( \normaHH{(\phieps,\phiGeps)}  + 1 \bigr)
  \non
  \\
  & \quad {}
  + \normaHH{(\phieps,\phiGeps)}^2
  + 1
  \bigr\}
  \quad \aet \,.
  \label{prelim}
\end{align}

At this point, we are ready to prove the basic estimates.

\step
First a priori estimate

We test \eqref{primaeps} by $(\mueps,\muGeps)$.
At the same time, we test \eqref{secondaeps} by $(\dt\phieps,\dt\phiGeps)+\gamma(\phieps,\phiGeps)$.
Then, we sum up and observe \anold{that} several \betti{cancellations} occur.
Hence, we have, \an{almost everywhere in $(0,T)$,}~that
\Bsist
  && \iO |\nabla\mueps|^2
  + \iG |\nablaG\muGeps|^2
  + \tau \iO |\dt\phieps|^2
  + \tau \iG |\dt\phiGeps|^2
  \non
  \\
  && \quad {}
  + \frac 12 \, \frac d{dt} \iO |\nabla\phieps|^2
  + \frac 12 \, \frac d{dt} \iG |\nablaG\phiGeps|^2
  + \frac d{dt} \iO \Betaeps(\phieps)
  + \frac d{dt} \iG \BetaGeps(\phiGeps)
  \non
  \\
  && \quad {}
  + \iO \pi(\phieps) \dt\phieps
  + \iG \piG(\phiGeps) \dt\phiGeps
  \non
  \\
  && \quad {}
  + \frac {\tau\gamma} 2 \, \frac d{dt} \iO |\phieps|^2
  + \frac {\tau\gamma} 2 \, \frac d{dt} \iG |\phiGeps|^2
  + \gamma \iO |\nabla\phieps|^2
  + \gamma \iG |\nabla\phiGeps|^2
  \non
  \\
  && \quad {}
  + \gamma \iO \betaeps(\phieps) \phieps
  + \gamma \iG \betaGeps(\phiGeps) \anold{\phiGeps}
  + \gamma \iO \pi(\phieps) \phieps
  + \gamma \iG \piG(\phiGeps) \anold{\phiGeps}
  \non
  \\
  \separa
  && = \gamma \iO u \mueps
  + \gamma \anold{\iG} \uG \muGeps \,.
  \non
\Esist
Two integrals are nonnegative since $\betaeps(r)$, $\betaGeps(r)$ and $r$ have the same sign for every $r\in\erre$.
By ignoring them, integrating over $(0,t)$, where $t\in(0,T)$ is arbitrary, 
accounting for the initial condition~\eqref{cauchyeps} and
the \Lip\ continuity of $\pi$ and $\piG$, 
and applying the Young inequality with an arbitrary $\delta>0$, 
we deduce~that
(see \eqref{defQtSt} for the notation)
\Bsist
  && \intQt |\nabla\mueps|^2
  + \intSt |\nablaG\muGeps|^2
  + \tau \intQt |\dt\phieps|^2
  + \tau \intSt |\dt\phiGeps|^2
  \non
  \\
  && \quad {}
  + \frac 12 \iO |\nabla\phieps(t)|^2
  + \frac 12 \iG |\nablaG\phiGeps(t)|^2
  + \iO \Betaeps(\phieps(t))
  + \iG \BetaGeps(\phiGeps(t))
  \non
  \\
  && \quad {}
  + \frac {\tau\gamma} 2 \iO |\phieps(t)|^2
  + \frac {\tau\gamma} 2 \iG |\phiGeps(t)|^2
  + \rev{\gamma} \intQt |\nabla\phieps|^2
  + \rev{\gamma} \intSt |\nablaG\phiGeps|^2
  \non
  \\
  \separa
  && \leq \frac 12 \iO |\nabla\phiz|^2
  + \frac 12 \iG |\nablaG\phiGz|^2
  + \iO \Betaeps(\phiz)
  + \iG \BetaGeps(\phiGz)
  \non
  \\
  && \quad {}
  + \frac {\tau\gamma} 2 \iO |\phiz|^2
  + \frac {\tau\gamma} 2 \iG |\phiGz|^2
  \non
  \\
  \separa
  && \quad {}
  + \delta \intQt |\dt\phieps|^2
  + \delta \intSt |\dt\phiGeps|^2
  + c_\delta \intQt (|\phieps|^2 + 1)
  + c_\delta \intSt (|\phiGeps|^2 + 1)
  \non
  \\
  && \quad {}
  + \gamma \intQt u \mueps
  + \gamma \intSt \uG \muGeps \,.
  \non
\Esist
By recalling the assumption \eqref{hpphiz} on the initial datum
and the properties of the Yosida approximation (see, e.g., \eqref{propBetaeps}),
it is clear that just the last two integrals need some treatment.
With the notation \eqref{defalpha} already used,
by also recalling \eqref{hpu} and the Poincar\'e inequality \eqref{poincare},
we have~that
\Bsist
  && \intQt u \mueps
  + \intSt \uG \muGeps
  \non
  \\
  && = \intQt u (\mueps-\aeps)
  + \intSt \uG (\muGeps-\aeps)
  + \iot \aeps(s) \Bigl( \iO u(s) + \iG \uG(s) \Bigr) \, ds
  \non
  \\
  && \leq \delta \intQt |\nabla\mueps|^2
  + \delta \intSt |\nablaG\muGeps|^2
  + c_\delta
  + c \iot |\aeps(s)| \, ds \,.
  \non
\Esist
On the other hand, by our preliminary estimate \eqref{prelim} and the Young and Poincar\'e inequalities once more, 
we have~that
\Bsist
  && \iot |\aeps(s)| \, ds
  \leq c \iot \bigl\{
    \normaHH{(\dt\phieps,\dt\phiGeps)} \, \bigl( \normaHH{(\phieps,\phiGeps)} + 1 \bigr)
  \non
  \\
  &&\quad {}
  + \normaHH{(\mueps-\aeps,\muGeps-\aeps)} \, \bigl( \normaHH{(\phieps,\phiGeps)}  + 1 \bigr)
  + \normaHH{(\phieps,\phiGeps)}^2
  + 1
  \bigr\} \, ds
  \non
  \\
  && \leq \delta \intQt |\dt\phieps|^2
  + \delta \intSt |\dt\phiGeps|^2
  + \delta \intQt |\nabla\an{\mueps}|^2
  + \delta \intSt |\nabla\an{\muGeps}|^2
  \non
  \\
  && \quad {}
  + c_\delta \intQt \bigl( |\phieps|^2 + 1 \bigr) 
  + c_\delta \intSt \bigl( |\phiGeps|^2 + 1 \bigr) .
  \non
\Esist
{Therefore}, choosing $\delta$ small enough and applying the Gronwall lemma,
{lead us to} conclude~that
\Beq
  \norma{(\phieps,\phiGeps)}_{\H1\calH\cap\L\infty\calV}
  + \norma{(\nabla\mueps,\nablaG\muGeps)}_{\L2\calH}
  \leq c \,.
  \label{primastima}
\Eeq

\step
Consequence

We come back to \eqref{prelim} and account for \eqref{primastima}.
We \anold{realize} that
\Bsist
  && \ioT |\aeps(s)|^2 \, ds
  \leq c \bigl\{
    \norma{(\dt\phieps,\dt\phiGeps)}_{\L2\calH}^2 \, \bigl( \norma{(\phieps,\phiGeps)}_{\L\infty\calH}^2 + 1 \bigr)
  \non
  \\
  && \quad {}
  + \norma{(\mueps-\aeps,\muGeps-\aeps)}_{\L2\calH}^2 \, \bigl( \norma{(\phieps,\phiGeps)}_{\L\infty\calH}^2  + 1 \bigr)
  \non
  \\
  && \quad {}
  + \norma{(\phieps,\phiGeps)}_{\L\infty\calH}^4
  + 1
  \bigr\} 
  \leq c \,.
  \non
\Esist
Therefore, \anold{using once more the Poincar\'e inequality,} we conclude that
\Beq
  \norma{(\mueps,\muGeps)}_{\L2\calV} \leq c \,.
  \label{stimamuL2V}
\Eeq

\step
Second a priori estimate

We test \eqref{secondaeps} by $(\betaeps(\phieps),\betaeps(\phiGeps))$
and integrate with respect to time.
After, rearranging, we obtain for every $t\in[0,T]$ that
\Bsist
  && \tau \iO \Betaeps(\phieps(t))
  + \tau \iG \Betaeps(\phiGeps(t))
  + \intQt \betaeps'(\phieps) |\nabla\phieps|^2
  + \intSt \betaeps'(\phiGeps) |\nablaG\phiGeps|^2
  \non
  \\
  && \quad {}
  + \intQt |\betaeps(\phieps)|^2
  + \intSt \betaGeps(\phiGeps) \, \betaeps(\phiGeps) 
  \non
  \\
  && = \tau \iO \Betaeps(\phiz)
  + \tau \iG \Betaeps(\phiGz)
  - \intQt \pi(\phieps) \, \betaeps(\phieps)
  \anold{-} \intSt \piG(\phiGeps) \, \betaeps(\phiGeps)
  \non
  \\
  && \quad {}
  + \anold{\intQt} \mueps \betaeps(\phieps)
  + \anold{\intSt} \muGeps \betaeps(\phiGeps) \,.
  \non
\Esist
{All of the terms on the \lhs\ are nonnegative.
However, as for the last one, we can say something better, namely,}
\Beq
  \intSt \betaGeps(\phiGeps) \, \betaeps(\phiGeps) 
  \geq \frac 1 {2C^*} \intSt |\betaeps(\phiGeps)|^2 - c
  \non
\Eeq
due to \eqref{prodbetaeps}.
This also makes easier to bound the \rhs,
on account of \eqref{propBetaeps},
\eqref{dapropCCeps},
\eqref{hpphiz},
the \Lip\ continuity of $\pi$ and~$\piG$,
the Young inequality and estimates \accorpa{primastima}{stimamuL2V}
\anold{ and allow us to} conclude~that
\Beq
  \norma{\betaeps(\phieps)}_{\L2H}
  + \norma{\betaeps(\phiGeps)}_{\L2\HG}
  \leq c \,.
  \label{secondastima}
\Eeq

\step
Third a priori estimate

We write \eqref{secondaeps} \aat\ in the form
\Bsist
  && \iO \nabla\phieps(t) \cdot \nabla v
  + \iG \nablaG\phiGeps(t) \cdot \nablaG\vG
  + \iG \betaGeps(\phiGeps(t)) \, \vG
  \non
  \\
  && = \iO g \, v
  + \iG \gG \, \vG
  \quad \hbox{for every $(v,\vG)\in\calV$}, 
  \quad \hbox{where}
  \non
  \\
  && g := \mueps(t) - \tau \dt\phieps(t) - (\betaeps + \pi)(\phieps(t))
  \aand
  \gG := \muGeps(t) - \tau \dt\phiGeps(t) - \piG(\phiGeps(t))
  \non
\Esist
and apply Lemma~\ref{Elliptic} with $\sigma=\betaGeps$.
As the constant $\CO$ appearing in \eqref{ellreg} does not depend on~$\eps$,
by squaring \eqref{ellreg} itself in our case, integrating over $(0,T)$ and accounting for the estimates already obtained,
we conclude~that
\Beq
  \norma{(\phieps,\phiGeps)}_{\L2\calW}
  + \norma{\betaGeps(\phiGeps)}_{\L2\HG}
  \leq c \,.
  \label{terzastima}  
\Eeq
Similarly, by applying the lemma to \eqref{primaeps}, we obtain that
\Beq
  \norma{(\mueps,\muGeps)}_{\L2\calW}
  \leq c \,.
  \label{terzastimabis}  
\Eeq

\step
Conclusion of the existence proof

By collecting the estimates we have established 
and applying \an{well-known} weak and weak star compactness results,
we deduce~that 
\Bsist
  && (\mueps,\muGeps) \to (\mu,\muG)
  \quad \hbox{weakly in $\L2\calW$}
  \non
  \\
  && (\phieps,\phiGeps) \to (\phi,\phiG)
  \quad \hbox{weakly star in $\H1\calH\cap\L\infty\calV\cap\L2\calW$}
  \non
  \\
  && (\betaeps(\phieps),\betaGeps(\phiGeps)) \to (\xi,\xiG)
  \quad \hbox{weakly in $\L2\calH$}
  \non
\Esist
(at~least for a subsequence $\eps_k\seto0$)
for some limit {six}-tuple $\soluz$.
Then it is clear that the integrated versions of \eqref{prima} and \eqref{seconda}
with time dependent test functions $(v,\vG)$ in $\L2\calV$ are satisfied.
Moreover, $\phi$~clearly verifies the initial condition \eqref{cauchy}.
Finally, the conditions \eqref{terza} hold as well 
since $(\phieps,\phiGeps)$ converges strongly, 
e.g., in $\L2\calH$ by the Aubin--Lions lemma
(see, e.g., \cite[Thm.~5.1, p.~58]{Lions})
and one can apply \an{well-known} results on monotone operators 
(see, e.g., \cite[Lemma 2.3, p.~38]{Barbu}) \anold{to infer that $\xi \in \beta (\phi)$ and  $\xiG \in \betaG (\phiG)$, as
desired}.
This concludes the proof of Theorem~\ref{Existence}.


\subsection{Regularity and separation}
\label{REG-SEP}

We first prove Theorem~\ref{Regularity}.
Due to the existence proof just given, 
it is clear that it suffices to establish \anold{further} a priori estimates on the solution to the approximating problem
that correspond\an{s} to the regularity \eqref{regularity} and the inequality~\eqref{stabilitybis}.
Also in this case, we proceed formally acting on the solution $\soluzeps$, directly.
We notice that the time derivative and the initial value of the approximating chemical potential do not exist.
However, they do exist for the discrete solution given by a Faedo--Galerkin scheme.
Indeed, the discrete phase field solves a systems of \anold{ODEs} and is \anold{in fact} much smoother,
whence the discrete chemical potential is smoother too, as a consequence of the discrete analogue of~\eqref{secondaeps}.
Thus, the estimates we prove just formally would be correct when performed on the discrete solution.
Of course, in realizing our project, we account for the estimates already established in the previous section.
As before, $\delta$~\an{indicates} a positive parameter whose value is chosen when it is necessary to do~\an{so}.

\step
Fourth a priori estimate

We test \eqref{primaeps} by $(\dt\mueps,\dt\muGeps)$.
At the same time, we differentiate \eqref{secondaeps} with respect to time
and test the equality we obtain by $(\dt\phieps,\dt\phiGeps)$.
Then, we sum up and observe that some cancellation occur \anold{obtaining that}
\Bsist
  && \frac 12 \, \frac d{dt} \iO |\nabla\mueps|^2
  + \frac 12 \, \frac d{dt} \iG |\nablaG\muGeps|^2
  + \frac 12 \, \frac d{dt} \iO |\dt\phieps|^2
  + \frac 12 \, \frac d{dt} \iG |\dt\phiGeps|^2
  \non
  \\
  && \quad {}
  + \iO |\nabla\dt\phieps|^2
  + \iG |\nablaG\dt\phiGeps|^2
  + \iO \betaeps'(\phieps) |\dt\phieps|^2
  + \iG \betaGeps'(\phiGeps) |\dt\phiGeps|^2
  \non
  \\
  && = \gamma \iO (u-\phieps) \dt\mueps
  + \gamma \iG (\uG-\phiGeps) \dt\muGeps
  \anold{-} \iO \pi'(\phieps) |\dt\phieps|^2
  \anold{-} \iG \piG'(\phiGeps) |\dt\phiGeps|^2 \,.
  \non
\Esist
At this point, we integrate over $(0,t)$ with an arbitrary $t\in[0,T]$.
The volume integrals involving $\anold{\pi'}$ and $\anold{\piG'}$ 
are bounded due to {the previous estimates} \anold{and to \eqref{hpPi}}.
On the contrary, the initial values of the chemical potentials 
and those of the time derivatives of the phase variables that appear
have to be estimated.
At the same time, we have to \an{control} the integrals 
\an{arising} from the terms involving $\dt\mueps$ and~$\dt\muGeps$.
Let us first consider the latter, {which can be handled by integration by parts
as done below, where we also subtract and add the same quantities for convenience}.
With the notation {in}~\eqref{defalpha}, we obtain~that
\Bsist
  && \intQt (u-\phieps) \dt\mueps
  + \intSt (\uG-\phiGeps) \dt\muGeps
  \non
  \\
  \separa
  && = - \intQt (\dt u - \dt\phieps) \mueps
  - \intSt (\dt\uG - \dt\phiGeps) \muGeps
  \non
  \\
  && \quad {}
  + \iO \bigl( u(t) - \phieps(t) \bigr) \bigl( \mueps(t) - \aeps(t) \bigr)
  + \iG \bigl( \uG(t) - \phiGeps(t) \bigr) \bigl( \muGeps(t) - \aeps(t) \bigr)
  \non
  \\
  && \quad {}
  + \aeps(t) \Bigl( \iO \bigl( u(t) - \phieps(t) \bigr) + \iG \bigl( \uG(t) - \phiGeps(t) \bigr) \Bigr)
  \non
  \\
  && \quad {}
  - \iO \bigl( u(0) - \phiz \bigr) \mueps(0)
  - \iO \bigl( \uG(0) - \phiGz \bigr) \muGeps(0).
  \non
\Esist
The first two integrals on the \rhs\ are bounded due to our assumption \eqref{hpureg} on $(u,\uG)$
and to the estimates \eqref{primastima} and \eqref{stimamuL2V} of the previous section.
The whole second line of the \rhs, which we term $\grI$ for brevity, 
can be treated by the Young and Poincar\'e inequalities 
and it is bounded from above as follows
\Beq
  \grI \leq \delta \norma{(\nabla\mueps(t),\nablaG\muGeps(t))}_\calH^2
  + c_\delta \,.
  \non
\Eeq
As for the next line, it suffices to estimate the mean value.
By recalling \eqref{prelim} and \eqref{primastima}, we have~that 
\Beq
  |\aeps(t)|
  \leq \delta \normaHH{(\dt\phieps(t),\dt\phiGeps(t))}^2 
  + \delta \normaHH{(\nabla\mueps(t),\nablaG\muGeps(t))}^2
  + c_\delta \,.
  \label{peraepsbdd}
\Eeq
Finally, the last line is under control once the initial value of $(\mu,\muG)$ is estimated in~$\calH$.
Now, we proceed in estimating all the initial values we have to consider.
We write \eqref{primaeps} and \eqref{secondaeps} at the time $t=0$ and have~that
\begin{align}
  & \iO \dt\phieps(0) \, v
  + \iG \dt\phiGeps(0) \, \vG
  + \iO \nabla\mueps(0) \cdot \nabla v
  + \iG \nablaG\muGeps(0) \cdot \nablaG\vG
  \non
  \\
  & = \iO f_\eps v
  + \iG \fG_\eps \vG
  \label{primaepsz}
  \\
  & \tau \iO \dt\phieps(0) \, v
  + \tau \iG \dt\phiGeps(0) \, \vG
  + \iO \nabla\phiz \cdot \nabla v
  + \iG \nablaG\phiGz \cdot \nablaG\vG
  \non
  \\
  & = \iO \mueps(0) v
  + \iG \muGeps(0) \vG
  + \iO g_\eps v 
  + \iG \gG_\eps \vG
  \label{secondaepsz}
\end{align}
both for every $(v,\vG)\in\calV$,
with $f_\eps$ and $g_\eps$ bounded in $H$ and $\fG_\eps$ and $\gG_\eps$ bounded in~$\HG$
(in~particular, thanks to \eqref{propbetaeps} and \eqref{hpphizreg}).
Now, we subtract \eqref{secondaepsz} from \eqref{primaepsz} multiplied by $\tau$ and deduce~that
\Bsist
  && \tau \iO \nabla\mueps(0) \cdot \nabla v
  + \tau \iG \nablaG\muGeps(0) \cdot \nablaG\vG
  + \iO \mueps(0) v
  + \iG \muGeps(0) \vG
  \non
  \\
  && = \iO (\tau f_\eps - g_\eps) v
  + \iG (\tau\fG_\eps -\gG_\eps) \vG
  \anold{+} \iO \nabla\phiz \cdot \nabla v
  \anold{+} \iG \nablaG\phiGz \cdot \nablaG\vG
  \non
\Esist
still for every $(v,\vG)\in\calV$.
By applying the Lax--Milgram lemma
(or the Riesz representation theorem with an equivalent inner product in~$\calV$),
we obtain the estimate
\Beq
  \normaVV{(\mueps(0),\muGeps(0))} \leq c \,.
  \non
\Eeq
Now, by coming back to \eqref{secondaepsz} and recalling that $(\phiz,\phiGz)\in\calW$, 
we infer that
\Bsist
  && \tau \iO \dt\phieps(0) \, v
  + \tau \iG \dt\phiGeps(0) \, \vG
  \non
  \\
  && = \iO \bigl( \mueps(0) + g_\eps + \Delta\phiz \bigr) v
  \rev{{} - \iG \dn \phiz \, \vG}
  + \iG \bigl( \muGeps(0) + \gG_\eps + \DeltaG(\phiGz) \bigr)\vG
  \non
\Esist
and deduce that
\Beq
  \normaHH{(\dt\phieps(0),\dt\phiGeps(0))} \leq c \,.
  \non
\Eeq
By collecting all this, choosing $\delta$ small enough and applying the Gronwall lemma, we conclude that
\Beq
  \norma{(\nabla\mueps,\nablaG\muGeps)}_{\L\infty\calH}
  + \norma{(\dt\phieps,\dt\phiGeps)}_{\L\infty\calH\cap\L2\calV}
  \leq c \,.
  \label{quartastima}
\Eeq
On the other hand, we also deduce from this estimate and \eqref {peraepsbdd} 
that $\aeps$ is bounded in $L^\infty(0,T)$.
Therefore, we improve the above estimate by obtaining
\Beq
  \norma{(\mueps,\muGeps)}_{\L\infty\calV}
  \leq c \,.
  \label{daquartastima}
\Eeq

\step
Fifth {a priori} estimate

\anold{Next, we argue as done in the derivation of \eqref{secondastima},} but avoid integration in time.
We have \aet~that
\Bsist
  && \iO \betaeps'(\phieps) |\nabla\phieps|^2
  + \iG \betaeps'(\phiGeps) |\nablaG\phiGeps|^2
  + \iO |\betaeps(\phieps)|^2 
  + \iG \betaGeps(\phiGeps) \, \betaeps(\phiGeps)
  \non
  \\
  && = \iO \bigl( \mueps - \tau \dt\phieps - \pi(\phieps) \bigl) \betaeps(\phieps)
  + \iG \bigl( \muGeps - \tau \dt\phiGeps - \piG(\phiGeps) \bigl) \betaeps(\phiGeps)
  \non
\Esist
and accounting for \eqref{hpCC} and \accorpa{quartastima}{daquartastima}\anold{,}
we easily deduce that
\Beq
  \norma{\betaeps(\phieps)}_{\L\infty H}
  + \norma{\betaeps(\phiGeps)}_{\L\infty\HG}
  \leq c \,.
  \label{quintastima}
\Eeq

\step
Sixth a priori estimate

At this point, we can read \eqref{primaeps} written \aat\
as a particular case of \eqref{elliptic} with $\sigma\an{\equiv}0$,
apply Lemma~\ref{Elliptic} and account for the $L^\infty$ estimates just obtained.
We conclude~that
\Beq
  \norma{(\phieps,\phiGeps)}_{\L\infty\calW}
  + \norma{\betaGeps(\phiGeps)}_{\L\infty\HG}
  \leq c \,.
  \label{sestastima}  
\Eeq
Similarly, by applying the lemma to \eqref{primaeps}, we obtain that
\Beq
  \norma{(\mueps,\muGeps)}_{\L\infty\calW}
  \leq c \,.
  \label{sestastimabis}  
\Eeq

\step
Conclusion of the proof of Theorem~\ref{Regularity}

By proceeding as at the end of the proof of Theorem~\ref{Existence} given in Section~\ref{EXISTENCE},
we see that the approximating solution converges 
(as~$\eps\seto0$ at least for a subsequence)
to~a solution to problem \Pbl\ also in the weak star topology 
associated to the regularity requirements \eqref{regularity}
and that this solution satisfies estimate \eqref{stabilitybis} with a constant $K_2$ 
with the dependence specified in the statement.
\QED

\medskip

Now, we \anold{move to} {prove} Theorem~\ref{Separation}.
Recalling \eqref{hpD}, it is clear that the thesis in the first case
trivially follows from the boundedness of the component $\phi$ ensured by Theorem~\ref{Regularity}.
Thus, we suppose that $D$ is a bounded open interval.
In fact, we prove that every solution $\sol$ satisfies the separation property
whenever the components $\mu$ and $\muG$ are bounded.
In particular, this holds for the solution given by Theorem~\ref{Regularity}.
We assume $D=(-1,1)$, without loss of generality,
and notice that, 
since $\beta$ and $\betaG$ are maximal monotone and $\pi$ 
and $\piG$ are smooth in the whole of~$\erre$,
we have~that
\Beq
  \lim_{r\seto -1^+} F'(r)
  = \lim_{r\seto -1^+} \FG'(r)
  = - \infty
  \aand
  \lim_{r\neto 1^-} F'(r)
  = \lim_{r\neto 1^-} \FG'(r)
  = + \infty \,.
  \label{dahpD}
\Eeq
Hence, thanks to assumption \eqref{hpsepar} on $\phiz$
and the boundedness of $(\mu,\muG)$ we are assuming,
we can fix a real number $\rz$ satisfying
\begin{align}
  & \norma\phiz_\infty \leq \rz < 1 
  \label{hprz}
  \\
  \noalign{\noindent and, with $N := \rev{ \norma\mu_{L^\infty(Q)}}\,,$}
  & F'(r) \leq  -N
  \aand
  \FG'(r) \leq  -N
  \quad \hbox{if $r\in(-1,-\rz]$}
  \label{menoinfty}
  \\
  & F'(r) \geq  N
  \aand
  \FG'(r) \geq  N
  \quad \hbox{if $r\in[\rz,1)$}.
  \label{piuinfty}
\end{align}
Then, we choose a function $\zeta:\erre\to\erre$ that is monotone and \Lip\ continuous and satisfies
\Beq
  \zeta(r) = 0 
  \quad \hbox{if $|r|\leq\rz$} , \quad
  \zeta(r) < 0 
  \quad \hbox{if $r<-\rz$}
  \aand
  \zeta(r) > 0
  \quad \hbox{if $r>\rz$}
  \non
\Eeq
and test \eqref{seconda} by $(\zeta(\phi),\zeta(\phiG))$.
By setting $\hat\zeta( r):=\int_0^{ r}\zeta(s)\,ds$ for $r\in\erre$, we obtain~that
\Bsist
  && \tau \, \frac d{dt} \Bigl(
    \iO \hat\zeta(\phi)
    + \iG \hat\zeta(\phiG)
  \Bigr)
  + \iO \zeta'(\phi) |\nabla\phi|^2
  + \iG \zeta'(\phiG) |\nablaG\phiG|^2
  \non
  \\
  && = \iO \bigl( \mu - F'(\phi) \bigr) \zeta(\phi)
  + \iG \bigl( \muG - \FG'(\phiG) \bigr) \zeta(\phiG)
  \leq 0 \,.
  \non
\Esist
We \anold{then} deduce that
\Beq
  \iO \hat\zeta(\phi(t))
  + \iG \hat\zeta(\phiG(t))
  \leq \iO \hat\zeta(\phiz)
  + \iG \hat\zeta(\phiGz)
  = 0 
  \non
\Eeq
for every $t\in[0,T]$.
Since $\hat\zeta$ is nonnegative, we conclude that $\zeta(\phi)$ vanishes identically,
i.e., that $|\phi|\leq\rz$ in the whole of~$Q$ \anold{leading to the existence of $\phi_*$ and $\phi^*$ as in the statement}.
\QED


\section{Continuous dependence and uniqueness}
\label{CONTDEP}
\setcounter{equation}{0}

This section is devoted to the proof of Theorem~\ref{Contdep}.
More precisely, we prove estimate \eqref{contdep}
for an arbitrary pair of solutions, so that uniqueness follows as a consequence.
\anold{Let us} fix two pairs $(u_i,u_{\Gamma,i})$, $i=1,2$, of \anold{forcing terms}
and consider arbitrary corresponding solutions $(\mu_i,\mu_{\Gamma,i},\phi_i,\phi_{\Gamma,i})$.
{For brevity,} we term $(u,\uG)$ and $\sol$ the differences.
Then, $\sol$ solves the equations
\begin{align}
  & \iO \dt\phi \, v
  + \iG \dt\phiG \, \vG
  + \iO \nabla\mu \cdot \nabla v
  + \iG \nablaG\muG \cdot \nablaG\vG
  \non
  \\
  & = \gamma \iO (u-\phi) v
  + \gamma \iG (\uG-\phiG) \vG
  \label{diffprima}
  \\
  \separa
  & \tau \iO \dt\phi \, v
  + \tau \iG \dt\phiG \, \vG
  + \iO \nabla\phi \cdot \nabla v
  + \iG \nablaG\phiG \cdot \nablaG\vG
  \non
  \\
  & \quad {}
  + \iO \bigl( F'(\phi_1)-F'(\phi_2) \bigr) v
  + \iG \bigl( \FG'(\phi_{\Gamma,1})-\FG'(\phi_{\Gamma,2}\anold{)} \bigr) \vG
  = \iO \mu v 
  + \iG \muG \vG
  \label{diffseconda}
\end{align}
both \aet\ and for every $(v,\vG)\in\calV$,
as well as the initial condition $\phi(0)=0$.
We test the above equations by $(\mu,\muG)$ and $(\dt\phi,\dt\phiG)+\gamma(\phi,\phiG)$, respectively.
Then, we sum up and notice that a number of cancellation occur.
Moreover, we recall that the values of $\phi_1$ and $\phi_2$ 
belong to the interval $[\phimin,\phimax]$ given by Theorem~\ref{Separation}
and that $F$ and $\FG$ are \Lip\ continuous in this interval.
Thus, we can use this fact in treating the terms involving $F$ and~$\FG$ moved to the \rhs.
By also owing to the Young inequality, we have~that
\begin{align}
  & \iO |\nabla\mu|^2
  + \iG |\nablaG\muG|^2
  \non
  \\
  \separa
  & \quad {}
  + \tau \iO |\dt\phi|^2
  + \tau \iG |\dt\phiG|^2
  + \frac 12 \, \frac d{dt} \iO |\nabla\phi|^2
  + \frac 12 \, \frac d{dt} \iG |\nablaG\phiG|^2
  \non
  \\
  & \quad {}
  + \frac {\tau\gamma} 2 \, \frac d{dt} \iO |\phi|^2  
  + \frac {\tau\gamma} 2 \, \frac d{dt} \iG |\phiG|^2
  + \gamma \iO |\nabla\phi|^2
  + \gamma \iG |\nablaG\phiG|^2
  \non
  \\
  & \leq \gamma \iO u \mu
  + \gamma \iG \uG \muG 
  + \frac \tau 4 \iO |\dt\phi|^2
  + \frac \tau 4 \iG |\dt\phiG|^2
  + c \iO |\phi|^2
  + c \iG |\phiG|^2 \,.
  \label{percontdep}
\end{align}
Now, we test \eqref{seconda} by $(1,1)$.
By setting $\alpha:=\mean(\mu,\muG)$, 
and using the \Lip\ continuity of $F$ and $\FG$ once more, we obtain
\Beq
  (\misO+\misG) |\alpha|
  \leq \tau \iO |\dt\phi|
  + \tau \iG |\rev{\dt\phiG}|
  + c \iO |\phi|
  + c \iG |\phiG|
  \non
\Eeq
whence
\Bsist
  && (\misO+\misG)^2 |\alpha|^2
  \leq 4\tau^2 \Bigl( \iO |\dt\phi| \Bigr)^2
  + 4\tau^2 \Bigl( \iG |\dt\phiG| \Bigr)^2
  + c \Bigl( \iO |\phi| \Bigr)^2
  + c \Bigl( \iG |\phiG| \Bigr)^2
  \non
  \\
  && \leq 4\tau^2\misO \iO |\dt\phi|^2
  + 4\tau^2\misG \iG |\dt\phiG|^2
  + c \iO |\phi|^2
  + c \iG |\phiG|^2 \,.
  \non
\Esist
By \anold{dividing} by $16\tau(\misO+\misG)$, we conclude that
\Beq  
  \anold{\frac {\misO+\misG} {16\tau}} \, |\alpha|^2
  \leq \frac \tau 4 \iO |\dt\phi|^2
  + \frac \tau 4 \iG |\dt\phiG|^2 
  + c \iO |\phi|^2
  + c \iG |\phiG|^2 \,.
  \non
\Eeq
Next, we add this inequality to \eqref{percontdep} and rearrange.
At this point, we see that we can treat the first two integrals on the \rhs\ 
by accounting for the Young and Poincar\'e inequalities \eqref{young} and~\eqref{poincare}.
For every $\delta>0$, we have indeed~that
\Bsist
  && \gamma \iO u \mu
  + \gamma \iG \uG \muG 
  \leq \delta \iO |\mu^2|
  + \delta \iG |\muG|^2
  + c_\delta \iO |u|^2
  + c_\delta \iG |\uG|^2 \,.
  \non
  \\
  && \leq \delta \CO \Bigl(
    \iO |\nabla\mu|^2
    + \iG |\nablaG\muG|^2
    + |\alpha|^2
  \Bigr)
  + c_\delta \iO |u|^2
  + c_\delta \iG |\uG|^2 \,.
  \non
\Esist
Therefore, by choosing $\delta$ small enough, integrating in time and applying the Gronwall lemma,
we obtain \eqref{contdep} with a constant $\rev{K_3}$ as in the statement.


\section{The control problem}
\label{CONTROL}
\setcounter{equation}{0}

In this section, we study the control problem \ControlPbl\ 
and prove the results we have stated in Section~\ref{STATEMENT}.
\anold{As already anticipated, the controls $u$ and $ \uG$ enter the state system \Pbl\ in the form of distributed and boundary source terms: of course, all the results stated above continue to hold for every admissible pair of controls, i.e., for every $(u,\uG) \in \Uad$.}
In the whole section, it is understood that all of the assumptions we have listed there are satisfied.
However, before starting, it is convenient to introduce a \nbh\ of the control box
and make an observation.
For $R>0$, we~set
\begin{align}
  & \UR := \bigl\{(u,\uG)\in\calU:
  \non
  \\
  & \phantom{\UR:=.\ } \norma u_\infty < M+R \ \ \hbox{and} \ \ \norma\uG_\infty < M+R,
  \non
  \\
  & \phantom{\UR:=.\ } \norma{\dt u}_{\L2H} < M'+R \ \ \hbox{and} \ \ \norma{\dt\uG}_{\L2\HG} < M'+R\bigr\}
  \label{defUR}
  \\
  & \hbox{where} \quad 
  \calU := \bigl(\LQ\infty\times\LS\infty\bigr)\cap\H1\calH
  \label{defU}
\end{align}
and we fix $R$ \anold{ small} in order that
\Beq
  \hbox{assumption \eqref{hpmz} is satisfied with $M$ replaced by $M+R$}.
  \label{hpR}
\Eeq
Therefore, all our results hold for every $(u,\uG)\in\UR$
and the constants of the stability estimates
(now depending on $R$ in addition)
can be fixed once and for all since $R$ is fixed.
In particular, the boundedness of the component $\phi$ in the case of regular potentials
or the separation property in the case of logarithmic type potential
combined with the regularity condition \eqref{hpF} ensure~that
\Beq
  \norma{F^{(i)}(\phi)}_\infty \leq K
  \aand
  \norma{\FG^{(i)}(\phiG)}_\infty \leq K
  \quad \hbox{for $0\leq i\leq3$}
  \label{start}
\Eeq
with a fixed constant~$K$,
where $\phi$ and $\phiG$ are the components of the solution to problem \Pbl\
corresponding to any $(u,\uG)\in\UR$.
For the same reason, we can assume~that 
\Beq
  \hbox{$F^{(i)}$ and $\FG^{(i)}$ are \Lip\ continuous for $0\leq i\leq2$}
  \label{potlip}
\Eeq
without loss of generality.


\subsection{Existence of an optimal \rev{control}}
\label{EXISTCONTROL}
 
We prove Theorem~\ref{OKcontrol} by using the direct method \anold{of calculus of variations}.
Thus, we pick a minimizing sequence $\graffe{(\un,\uGn)}$ in $\Uad$
and the corresponding sequence of the solutions,
labeled with the \rev{subscripts}~$n$, to the state system.
Since $\Uad$ is bounded in~$\calU$, we can assume~that
\Beq
  (\un,\uGn) \to (\ustar,\uGstar)
  \quad \hbox{weakly star in $\calU$}
  \non
\Eeq
as $n\neto\infty$ for some limit pair $(\ustar,\uGstar)$, which must belong to~$\Uad$,
{since it is convex and closed}.
{Besides}, the corresponding solutions are bounded as well 
in the topologies associated \an{with} our well-posedness and regularity results.
Hence, {it follows}~that
\Bsist
  && (\phin,\phiGn) \to (\phi,\phiG)
  \quad \hbox{weakly star in $\W{1,\infty}\calH\cap\H1\calV\cap\L\infty\calW$}
  \non
  \\
  && (\mun,\muGn) \to (\mu,\muG)
  \quad \hbox{weakly star in $\L\infty\calW$}
  \non
\Esist
for some limit functions.
By also accounting for \eqref{potlip}, it is immediately seen that the quadruplet
$\sol$ is the solution to the state system
corresponding to $(\ustar,\uGstar)$ and~that
\Beq
  \calJ(\un,\uGn;\phin,\phiGn) \to \calJ(\ustar,\uGstar;\phi,\phiG) \,.
  \non
\Eeq
{On the other hand, we also have that
$\lim_{n\neto\infty} \calJ(\un,\uGn;\phin,\phiGn)$
coincides with the infimum of~$\calJ$
since the sequence $\graffe{(\un,\uGn)}$ is minimizing~$\calJ$.
Therefore, $(\ustar,\uGstar)$ is an optimal control.}


\subsection{The control-to-state mapping}
\label{CTS}

In this section, we introduce the {\it control-to-state} mapping and prove its \frechet\ differentiability.
\an{Along with} the space $\calU$ defined by~\eqref{defU}, we introduce {the state space}
\Beq
  \calY := \H1\calH \cap \L\infty\calV 
  \label{defY}
\Eeq
and the \anold{solution mapping} $\calS:\UR\to\calY$ defined~by
\Beq
  \vbox{\hsize.8\hsize\noindent
  for $(u,\uG)\in\UR$, the equality $(\phi,\phiG)=\calS(u,\uG)$ means that
  $\phi$ and $\phiG$ are the components of the solution $\sol$ to \Pbl\ corresponding to $(u,\uG)$.}
  \label{defS}
\Eeq
Thanks to all the observations we have made before, this map is well-defined.
The main result of this section is \anold{showing} its \frechet\ \betti{differentiability}.
This is related to the linearized system we introduce at once.
Let $(u,\uG)\in\UR$ and let $\sol$ be the corresponding solution to problem \Pbl.
Then, the linearized system corresponding to $(u,\uG)$ and to the variation $(h,\hG)\in\calU$
is the system
\begin{align}
  & \iO \dt\psi \, v
  + \iG \dt\psiG \, \vG
  + \iO \nabla\eta \cdot \nabla v
  + \iG \nablaG\etaG \cdot \nablaG\vG
  \non
  \\
  & = \gamma\iO (h-\psi) v
  + \gamma \iG (\hG-\psiG) \vG
  \quad \hbox{\aet\ and for every $(v,\vG)\in\calV$}
  \qquad
  \label{primal}
  \\
  \separa
  & \tau \iO \dt\psi \, v
  + \tau \iG \dt\psiG \, \vG
  + \iO \nabla\psi \cdot \nabla v
  + \iG \nablaG\psiG \cdot \nablaG\vG
  \non
  \\
  & \quad {}
  + \iO F''(\phi) \psi v
  + \iG \FG''(\phiG) \psiG \vG
  \non
  \\
  & = \iO \eta v
  + \iG \etaG \vG
  \quad \hbox{\aet\ and for every $(v,\vG)\in\calV$}
  \label{secondal}
  \\
  \separa
  & \psi(0) = 0 
  \label{cauchyl}  
\end{align}
\Accorpa\Pbll primal cauchyl
where the unknown is the $\calV\times\calV$ valued function $\soluzl$.

One can show that this system has a unique solution satisfying
\Beq
  (\eta,\etaG) \in \L2\calV
  \aand
  (\psi,\psiG) \in \calY \,.
  \non
\Eeq
We do not provide the proof for brevity.
The existence of a solution can be proved by starting from a Faedo--Galerkin scheme.
{As the system \Pbll\ is linear},
uniqueness can be deduced as a trivial consequence of the stability estimate we establish at once.
Indeed, if $(h,\hG)=(0,0)$, the estimate implies that $(\eta,\etaG,\psi,\psiG)=(0,0,0,0)$.

\Blem
\label{PreF}
Let $(u,\uG)\in\UR$ and $(h,\hG)\in\calU$,
and let $\soluzl$ be any solution to the corresponding linearized system.
Then the estimate
\Beq
  {\norma{(\eta,\etaG)}_{\L2\calV}}
  + \norma{(\psi,\psiG)}_\calY
  \leq K \, \norma{(h,\hG)}_{\L2\calH}
  \label{preF}
\Eeq
holds true with a \anold{positive} constant $K$ that does not depend on~$(h,\hG)$.
\Elem

\Bdim
We just sketch the proof \anold{for brevity}.
We test equations \eqref{primal} and \eqref{secondal} by $(\eta,\etaG)$ 
and {$(\dt\psi,\dt\psiG)+\gamma(\psi,\psiG)-(1/\tau)(\eta,\etaG)$}, respectively,
and sum~up.
Then, several cancellation\an{s} occur.
Moreover, we move some terms from one side to the other one \anold{in~order to have just nonnegative integrals on the \lhs}
and repeatedly owe to the Young inequality on the \rhs\
by recalling that $F''(\phi)$ and $\FG''(\phiG)$ are bounded \anold{as a consequence of Theorem \ref{Separation} (or already of Theorem \ref{Regularity} in the case of everywhere defined potentials).}
Then, the only \anold{delicate term is handled as follows}
\Beq
  \iO \dt\rev{\psi }\, \eta
  \leq \frac \tau 2 \iO |\dt\rev{\psi}|^2
  + \frac 1 {2\tau} \iO |\eta|^2
  \non
\Eeq
and \an{similar} treatment is reserved to the analogous boundary \anold{contribution}.
Therefore, by integrating over $(0,t)$ with an arbitrary $t\in(0,T]$
and applying the Gronwall lemma,
we immediately obtain~\eqref{preF}.
\Edim

Here is {our} differentiability result.

\Bthm
\label{Frechet}
Given any $(u,\uG)\in\UR$,
the map $\calS$ is \frechet\ differentiable at $(u,\uG)$ 
and its \frechet\ derivative $D\calS(u,\uG)$ is the map belonging to $\calL(\calU,\calY)$ 
that to every $(h,\hG)\in\calU$ associates the pair $(\psi,\psiG)$
where $\soluzl$ is the solution to the linearized system
corresponding to $(u,\uG)$ and to the variation $(h,\hG)$.
\Ethm

\Bdim
The fact that the {linear} map described in the statement belongs to $\calL(\calU,\calY)$ 
is a consequence of Lemma~\ref{PreF}.
For the remaining part of the proof, without loss of generality,
we assume that $\norma{(h,\hG)}_\calU$ is small enough,
namely, {such that} the perturbed control pair $(u+h,\uG+\hG)$ also belongs to {the open set}~$\UR$.
We consider the solutions to the state system corresponding to $(u,\uG)$ and to $(u+h,\uG+\hG)$.
We term $\sol$ the former and $\solh$ the latter,
and define 
\Beq
  \theta := \muh - \mu - \eta, \quad
  \thetaG := \muGh - \muG - \etaG, \quad
  \rho := \phih - \phi - \psi
  \aand
  \rhoG := \phiGh - \phiG - \psiG \,.
  \non
\Eeq
Then, the quadruplet $(\theta,\thetaG,\rho,\rhoG)$ solves the system
\begin{align}
  & \iO \dt\rho \, v
  + \iG \dt\rhoG \, \vG
  + \iO \nabla\theta \cdot \nabla v
  + \iG \nablaG\thetaG \cdot \nablaG\vG
  \non
  \\
  & = - \gamma \iO \rho v
  - \gamma \iG \rhoG \vG
  \quad \hbox{\aet\ and for every $(v,\vG)\in\calV$}
  \qquad
  \label{primadiff}
  \\
  & \tau \iO \dt\rho \, v
  + \tau \iG \dt\rhoG \, \vG
  + \iO \nabla\rho \cdot \nabla v
  + \iG \nablaG\rhoG \cdot \nablaG\vG
  + \iO \chi v
  + \iG \chiG \vG
  \non
  \\
  & = \iO \theta v
  + \iG \thetaG \vG
  \quad \hbox{\aet\ and for every $(v,\vG)\in\calV$}
  \label{secondadiff}
  \\
  & \rho(0) = 0 
  \label{cauchydiff}  
\end{align}
where {we have set}
\Beq
  \chi := F'(\phih) - F'(\phi) - F''(\phi)\psi
  \aand
  \chiG := \FG'(\phiGh) - \FG'(\phiG) - \FG''(\phiG)\psiG \,.
  \label{defchi}
\Eeq
We notice at once that the Taylor formula with integral remainder yields that
\Bsist
  && \chi
  = F''(\phi)\rho + (\phih-\phi)^2 \, \calR
  \quad \hbox{where}
  \non
  \\
  &&  |\calR|
  = \Bigl| \int_0^1 (1-s) F\an{^{(3)}}(\phi+s(\phih-\phi)) \, ds \Bigr|
  \leq c \quad \aeQ 
  \non
\Esist
and that a similar observation holds for the boundary term~$\chiG$.
Therefore, {since $F''(\phi)$ and $\FG''(\phiG)$ are bounded},
for every $t\in(0,T]$, we have~that
\Beq
  \intQt |\chi|^2 + \intSt |\chiG|^2
  \leq c \intQt |\rho|^2 + c \intSt |\rhoG|^2
  + c \anold{\intQt} |\phih-\phi|^4 + c \anold{\intSt} |\phiGh-\phiG|^4 \,.
  \non
\Eeq
By accounting for inequality \eqref{contdep}, the continuous embedding
$\L\infty V\subset\LQ4$ and the boundary analogue, we deduce~that
\Beq
  \intQt |\chi|^2 + \intSt |\chiG|^2
  \leq c \intQt |\rho|^2 + c \intSt |\rhoG|^2
  + c \, \norma{(h,\hG)}_{\L2\calH}^4 \,.
  \label{stimachi}
\Eeq
At this point, we test \eqref{primadiff} by $(\theta,\thetaG)$ 
and \eqref{secondadiff} by $(\dt\rho,\dt\rhoG)+\gamma(\rho,\rhoG)$.
Then, we {add the resulting equalities to each other and notice that} 
some cancellations occur.
We obtain~that
\begin{align}
  & \iO |\nabla\theta|^2
  + \iG |\nablaG\thetaG|^2
  \non
  \\
  & \quad {}
  + \tau \iO |\dt\rho|^2
  + \tau \iG |\dt\rhoG|^2
  + \frac 12 \frac d{dt} \iO |\nabla\rho|^2
  + \frac 12 \frac d{dt} \iG |\nablaG\rhoG|^2
  \non
  \\
  & \quad {}
  + \rev{ \frac {\tau\gamma} 2} \, \frac d{dt} \iO |\rho|^2
  +\rev{ \frac {\tau\gamma} 2} \, \frac d{dt} \iG |\rhoG|^2
  + \gamma \iO |\nabla\rho|^2
  + \gamma \iG |\nablaG\rhoG|^2
  \non
  \\
  & = - \iO \chi \bigl( \dt\rho + \gamma\rho \bigl)
  - \iG \chiG \bigl( \dt\rhoG + \gamma\rhoG \bigl) \,.
  \label{perfrechet}
\end{align}
We aim at applying the Gronwall lemma.
To this end, we estimate the terms that \anold{arose} from the \rhs\
by integrating over~$(0,t)$.
They can be easily treated 
by accounting for the Young inequality and {the preliminary} estimate~\eqref{stimachi}.
\anold{Namely, we find that }
\Bsist
  && - \intQt \chi \bigl( \dt\rho + \gamma\rho \bigl)
  - \intSt \chiG \bigl( \dt\rhoG + \gamma\rhoG \bigl) 
  \non
  \\
  && \leq \frac \tau 2 \intQt |\dt\rho|^2
  + \frac \tau 2 \intSt |\dt\rhoG|^2
  + \intQt |\rho|^2
  + \intSt |\rhoG|^2
  + c \, \norma{(h,\hG)}_{\L2\calH}^4 \,.
  \non
\Esist
Therefore, we can rearrange \eqref{perfrechet} integrated over $(0,t)$ and apply the Gronwall lemma.
{This yields}
\Beq
  \norma{(\nabla\theta,\nablaG\thetaG)}_{\L2\calH}^2
  + \norma{(\rho,\rhoG)}_{\H1\calH\cap\L\infty\calV}^2
  \leq c \, \norma{(h,\hG)}_{\L2\calH}^4 
  \label{permoregen}
\Eeq
for every $(h,\hG)$ whose $\calU$-norm is small enough.
As a consequence, we have that
\Beq
  \norma{\calS(u+h,\uG+\hG)-\calS(u,\uG)-(\psi,\psiG)}
  = \norma{(\rho,\rhoG)}_\calY
  = o(\norma{(h,\hG)}_\calU \anold{)}
  \quad \hbox{as $\norma{(h,\hG)}_\calU\seto0$}
  \non
\Eeq
and this implies that $\calS$ is \frechet\ differentiable at $(u,\uG)$
and that its \frechet\ derivative acts as {expressed} in the statement.
\Edim

\Brem
\label{Moregeneral}
We observe that estimating $(\theta,\thetaG)$ from \eqref{secondadiff} with the help of~\eqref{permoregen}
yields~that
\Beq
  \norma{(\theta,\thetaG)}_{\L2\calH}^2
  \leq c \, \norma{(h,\hG)}_{\L2\calH}^4 \,.
  \non
\Eeq
This combined with \eqref{permoregen} itself and a corresponding improvement of \eqref{preF}
shows that also the map that to every $(u,\uG)\in\UR$ associates 
the whole solution $\sol$ to the original problem
is \frechet\ differentiable and that its \frechet\ derivative at the given $(u,\uG)$
is the map that to every $(h,\hG)\in\calU$
associates the whole solution $\soluzl$ to the linearized problem.
This fact would \anold{allow us \an{to consider}} a more general cost functional.
Namely, on the \rhs\ of \eqref{cost} we could add two similar integrals on $Q$ and $\Sigma$
involving $\mu$ and~$\muG$, respectively.
\Erem


\subsection{\anold{First-order} optimality conditions}
\label{OPT}

The above \frechet\ differentiability result
\anold{permits} us to apply the chain rule to the composite map
\Beq
  \UR \ni (u,\uG) \mapsto (u,\uG;\phi,\phiG) \mapsto \calJ(u,\uG;\phi,\phiG)
  \quad \hbox{with} \quad
  (\phi,\phiG) = \calS(u,\uG).
  \non
\Eeq
Since $\Uad$ is convex,
one immediately sees that a necessary condition for 
$(\ustar,\uGstar)$ to be an optimal control is given by the \anold{variational} inequality
\begin{align}
  & \alpha_1 \intQ (\phistar-\tarQ) \psi
  + \alpha_2 \intS (\phiGstar-\tarS) \psiG
  \non
  \\
  & \quad {}
  + \alpha_3 \iO (\phistar(T)-\tarO) \psi(T)
  + \alpha_4 \iG (\phiGstar(T)-\tarG) \psiG(T)
  \non
  \\
  & \quad {}
  + \alpha_5 \intQ \ustar (u-\ustar)
  + \alpha_6 \intS \uGstar (\uG-\uGstar)
  \non
  \\
  & \geq 0
  \quad \hbox{for every $(u,\uG)\in\Uad$}
  \label{preopt}
\end{align}
where $(\phistar,\phiGstar)=\calS(\ustar,\uGstar)$ and 
$\psi$ and $\psiG$ are the components of the solution $\soluzl$ to the linearized problem
corresponding to the pair $(\ustar,\uGstar)$ and 
to the variation $(h,\hG):=(u-\ustar,\uG-\uGstar)$.
However, the condition just written is \anold{rather} unpleasant.
Indeed, it involves the linearized problem infinitely many times,
since {$(u,\uG)$} is arbitrary in~$\Uad$.
As usual, this {issue} is bypassed by introducing a proper adjoint problem.

We give a rather weak formulation of it
since we do not require any time \an{differentiability} of the single components of its solution.
Namely, the only time derivative we consider is understood in the sense of the Hilbert triplet $(\calV,\calH,\calVp)$
obtained by identifying $\calH$ to a subspace of $\calVp$ in the usual way, i.e.,
in order that
$\<y,z>_\calV=(y,z)_\calH$ for every $y\in\calH$ and $z\in\calV$,
where $(\cpto,\cpto)_\calH$ denotes the standard inner product \anold{of}~$\calH$.
To make the presentation more readable, 
we introduce some abbreviations \anold{setting}
\Bsist
  && \lam := F''(\phistar) , \quad
  \lamG := \FG''(\phiGstar)
  \non
  \\
  && \zeta_1 := \alpha_1(\phistar-\tarQ) , \quad
  \zeta_2 := \alpha_2(\phiGstar-\tarS) 
  \non
  \\
  && \zeta_3 := \alpha_3(\phistar(T)-\tarO) , \quad
  \zeta_4 := \alpha_4(\phiGstar(T)-\tarG)
  \non
\Esist
where $(\ustar,\uGstar)$ is the optimal control at hand and $(\phistar,\phiGstar):=\calS(\ustar,\uGstar)$.
We recall that $\lam$ and $\lamG$ are bounded and that the other functions are \an{$L^2$} functions.
Then, the adjoint problem associated to $(\ustar,\uGstar)$ consists in looking for a quadruplet
(or~a pair of pairs) $(p,\pG,q,\qG)$ satisfying
\begin{align}
  & (p,\pG) \in \L\infty\calV
  \aand
  (q,\qG) \in \L\infty\calH \cap \L2\calV
  \label{regpq}
  \\
  & (p+\tau q,\pG+\tau\qG) \in \H1\calVp
  \label{dtptq}
\end{align}
\Accorpa\Regsoluza regpq dtptq
and solving the system
\begin{align}
  & - \< \dt(p+\tau q,\pG+\tau\qG) , (v,\vG) >_\calV
  \non
  \\
  & \quad {}
  + \iO \nabla q \cdot \nabla v
  + \iG \nablaG\qG \cdot \nablaG\vG
  + \iO (\gamma p + \lam q) v
  + \iG (\gamma\pG + \lamG\qG) \vG
  \non
  \\
  & = \iO \zeta_1 v
  + \iG \zeta_2 \vG
  \quad \hbox{\aet\ and for every $(v,\vG)\in\calV$}
  \label{primaa}
  \\
  & \iO \nabla p \cdot \nabla v
  + \iG \nablaG\pG \cdot \nablaG\vG
  \non
  \\
  & = \iO q v
  + \iG \qG \vG
  \quad \hbox{\aet\ and for every $(v,\vG)\in\calV$}
  \label{secondaa}
  \\
  & \rev{ \big( (p+\tau q,\pG+\tau\qG)(T) , (v,\vG) \big)_{\cal H}}
  \non
  \\
  & = \iO \zeta_3 v
  + \iG \zeta_4 \vG
  \quad \hbox{for every $(v,\vG)\in\calV$} .
  \label{cauchya}
\end{align}
\Accorpa\Pbla primaa cauchya

\Brem
We notice that \Regsoluza\ imply that
\Beq
  (p+\tau q,\pG+\tau\qG)
  \in \H1\calVp \cap \L2\calV
  \subset \C0\calH
  \non
\Eeq
so that the final value $(p+\tau q,\pG+\tau\qG)(T)$ is meaningful.
We also remark that the final condition \eqref{cauchya} is equivalent~to
\Beq
  (p+\tau q)(T) = \zeta_3
  \aand
  (\pG+\tau\qG)(T) = \zeta_4
  \non
\Eeq 
and that our assumptions on the cost functional simply yield that $(\zeta_3,\zeta_4)\in\calH$.
In particular, it is not required that $\zeta_4$ is the trace of~$\zeta_3$.
Such a condition would make sense only under stronger assumptions on the ingredients of the cost functional
and would be necessary if we were pretending $(p+\tau q,\pG+\tau\qG)$ to be a continuous $\calV$-valued function.
\Erem

At this point, we are ready to prove our results regarding the adjoint problem
and its connection with the control problem.

\Bthm
\label{Wellposednessadj}
The adjoint problem \Pbla\ has a unique solution $\soluza$ satisfying the regularity requirements \Regsoluza.
\Ethm

\Bdim
Our procedure is inspired by the proof of the \anold{comparable} result \cite[Thm.~4.4]{CGS14}.
As for existence,
we introduce an approximating problem depending on a positive parameter $\eps$ 
that we then let tend to zero.
This \anold{consists} in looking for a quadruplet $\soluzaeps$ satisfying
\Beq
  (\peps,\pGeps), \, (\qeps,\qGeps) \in \H1\calH \cap \L\infty\calV
  \label{regsoluzaeps}
\Eeq
and solving the system
\begin{align}
  & - \iO \dt(\peps+\tau\qeps) \, v
  - \iG \dt(\pGeps+\tau\qGeps) \, \vG
  \non
  \\
  & \quad {}
  + \iO \nabla\qeps \cdot \nabla v
  + \iG \nablaG\qGeps \cdot \nablaG\vG
  + \iO (\gamma\peps +\lam\qeps) v
  + \iG (\gamma\pGeps + \lamG\qGeps) \vG
  \non
  \\
  & = \iO \zeta_1 v
  + \iG \zeta_2 \vG
  \quad \hbox{\aet\ and for every $(v,\vG)\in\calV$}
  \label{primaaeps}
  \\
  & - \eps \iO \dt\peps \, v
  - \eps \iG \dt\pGeps \, \vG
  + \iO \nabla\peps \cdot \nabla v
  + \iG \nablaG\pGeps \cdot \nablaG\vG
  \non
  \\
  & = \iO \qeps v
  + \iG \qGeps \vG
  \quad \hbox{\aet\ and for every $(v,\vG)\in\calV$}
  \label{secondaaeps}
  \\
  & (\peps,\pGeps)(T) = (0,0)
  \aand
  (\qeps,\qGeps)(T)
  = (\zeta_3^\eps/\tau,\zeta_4^\eps/\tau)
  \label{cauchyaeps}
\end{align}
\Accorpa\Pblaeps primaaeps cauchyaeps
where the approximating data $\zeta_3^\eps$ and $\zeta_4^\eps$ are chosen such~that
\Beq
  (\zeta_3^\eps,\zeta_4^\eps) \in \calV
  \quad \hbox{for $\eps>0$}
  \aand
  (\zeta_3^\eps,\zeta_4^\eps) \to (\zeta_3,\zeta_4)
  \ \hbox{in $\calH$\quad as $\eps\seto0$}.
  \label{hpzetaeps}
\Eeq

\anold{Let us claim that this} problem has a unique solution.
However, we do not give the proof of this statement since it can be obtained just by closely following
the analogous one of \cite[Thm.~4.3]{CGS14},
where an ad~hoc inner product in $\calH$ is designed in order to make the problem parabolic.
On the contrary, 
we perform in detail the estimates that allow us to let $\eps$ tend to zero.
We assume $\eps\in(0,1)$ (so~that $\eps^2<\eps$) and employ the abbreviations
\Beq
  Q^t := \Omega\times(t,T)
  \aand
  \Sigma^t := \Gamma\times(t,T)
  \quad \hbox{for $t\in[0,T)$}.
  \non
\Eeq
\def\intQt{\int_{Q^t}}%
\def\intSt{\int_{\Sigma^t}}%
\def\intT{\int_t^T}%

\step
First a priori estimate

\anold{We test} \eqref{primaaeps} by $(\qeps,\qGeps)$ and \anold{rearrange the terms to find that}
\Bsist
  && - \iO \dt\peps \, \qeps
  - \iG \dt\pGeps \, \qGeps
  - \frac \tau 2 \, \frac d{dt} \iO |\qeps|^2
  - \frac \tau 2 \, \frac d{dt} \iG |\qGeps|^2
  \non
  \\
  && \quad {}
  + \iO |\nabla\qeps|^2
  + \iG |\nablaG\qGeps|^2
  \non
  \\
  && = - \iO (\gamma\peps +\lam\qeps) \qeps
  - \iG (\gamma\pGeps + \lamG\qGeps) \qGeps
  + \iO \zeta_1 \qeps
  + \iG \zeta_2 \qGeps \,.
  \non
\Esist
At the same time, by noting that $-\dt(\peps,\pGeps)$ belongs to $\L2\calV$
(since $(\qeps,\qGeps)\in\L2\calV$), we use it as test function in \eqref{secondaaeps} and obtain~that
\Bsist
  && \eps \iO |\dt\peps|^2
  + \eps \iG |\dt\pGeps|^2
  - \frac 12 \, \frac d{dt} \iO |\nabla\peps|^2
  - \frac 12 \, \frac d{dt} \iG |\nablaG\pGeps|^2
  \non
  \\
  && = - \iO \qeps \dt\peps
  - \iG \qGeps \dt\pGeps \,.
  \non
\Esist
Now, we sum up, account for the cancellations that occur,
and estimate the \rhs. \anold{Using that 
$\lam$ and $\lamG$ are bounded and that $\zeta_1$ and $\zeta_2$ are $L^2$ functions, we} deduce~that
\Bsist
  && - \frac \tau 2 \, \frac d{dt} \iO |\qeps|^2
  - \frac \tau 2 \, \frac d{dt} \iG |\qGeps|^2
  + \iO |\nabla\qeps|^2
  + \iG |\nablaG\qGeps|^2
  \non
  \\
  && \quad {}
  +  \eps \iO |\dt\peps|^2
  + \eps \iG |\dt\pGeps|^2
  - \frac 12 \, \frac d{dt} \iO |\nabla\peps|^2
  - \frac 12 \, \frac d{dt} \iG |\nablaG\pGeps|^2
  \non
  \\
  && \leq \iO |\peps|^2
  + \iG |\pGeps|^2
  + c \iO |\qeps|^2
  + c \iG |\qGeps|^2
  + c \,.
  \non
\Esist
At this point, we integrate over $(t,T)$.
As for the final values, we account for \eqref{cauchyaeps} and~\eqref{hpzetaeps}.
We obtain~that
\Bsist
  && \frac \tau 2 \iO |\qeps(t)|^2
  + \frac \tau 2 \iG |\qGeps(t)|^2
  + \intQt |\nabla\qeps|^2
  + \intSt |\nablaG\qGeps|^2
  \non
  \\
  && \quad {}
  + \eps \intQt |\dt\peps|^2
  + \eps \intSt |\dt\pGeps|^2
  + \frac 12 \iO |\nabla\peps(t)|^2
  + \frac 12 \iG |\nablaG\pGeps(t)|^2
  \non
  \\
  && \leq \intQt |\peps|^2 
  + \intSt |\pGeps|^2 
  + c \intQt |\qeps|^2 
  + c \intSt |\qGeps|^2 
  + c \,.
  \label{perprimastimaa}
\Esist
Now, we recall the Poincar\'e inequality \eqref{poincare} and have~that
\begin{align}
  & \intQt |\peps|^2 
  + \intSt |\pGeps|^2
  \non
  \\
  & \leq c \Bigl(
    \intQt |\nabla\peps|^2
    + \intSt |\nablaG\pGeps|^2
    + \intT |\mean(\peps,\pGeps)(s)|^2 \, ds
  \Bigr) .
  \label{dapoincare}
\end{align}
Next, we write the term $\gamma\peps$ in \eqref{primaaeps} as
$\gamma(\peps+\tau\qeps)-\gamma\tau\qeps$ 
and \rev{we do the same} with~$\gamma\pGeps$.
Moreover, we set for brevity
\Bsist
  && m := (\misO+\misG) \mean(\peps+\tau\qeps,\pGeps+\tau\qGeps)
  = \iO (\peps+\tau\qeps) + \iG (\pGeps+\tau\qGeps)
  \non
  \\
  && g := \zeta_1 + (\gamma\tau-\lam) \qeps
  \aand
  \gG := \zeta_2 + (\gamma\tau-\lamG) \qGeps \,.
  \non
\Esist
Then, we test both \eqref{primaaeps} with the above modification and \eqref{cauchyaeps}  by $(1,1)$ to obtain~that
\Beq
   - m' + \gamma m = \iO g + \iG \gG 
   \quad \aet
   \aand
   m(T) = \iO \zeta_3^\eps + \iG \zeta_4^\eps 
  \non
\Eeq
whence
\Beq
  m(t)
  = e^{-\gamma(T-t)} \Bigl( \iO \zeta_3^\eps + \iG \zeta_4^\eps \Bigr)
  + \intT e^{-\gamma(s-t)} \Bigl( \iO g(s) + \iG \gG(s) \Bigr) \, ds 
  \non
\Eeq
for every $t\in[0,T]$.
Therefore, by also recalling \eqref{hpzetaeps}, we easily conclude that
\Beq
  |\mean(\peps+\tau\qeps,\pGeps+\tau\qGeps)(t)|^2
  \leq c \intQt |\qeps|^2
  + c \intSt |\qGeps|^2 + c 
  \quad \hbox{for every $t\in[0,T]$}.
  \non
\Eeq
From this, we derive that
\begin{align}
  & |\mean(\peps,\pGeps)(t)|^2
  \leq C \, \Bigl( \iO |\qeps(t)|^2
  + \iG |\qGeps(t)|^2 \Bigr)
  \non
  \\
  & \quad {}
  + c \intSt |\qeps|^2
  + c \intSt |\qGeps|^2 + c 
  \quad \hbox{for every $t\in[0,T]$}
  \label{meanpeps}
\end{align}
\anold{with a fixed an computable positive constant $C$ depending} only on $\Omega$ and~$\tau$.
At this point, we add \eqref{perprimastimaa} and \eqref{meanpeps} multiplied by $\tau/(4C)$ to each other,
account for \eqref{dapoincare}, rearrange, and apply the (backward) Gronwall lemma.
By using the Poincar\'e inequality once more, we conclude~that
\Beq
  \norma{(\peps,\pGeps)}_{\L\infty\calV}
  + \norma{(\qeps,\qGeps)}_{\L\infty\calH\cap\L2\calV}
  + \eps^{1/2} \, \norma{\dt(\peps,\pGeps)}_{\L2\calH}
  \leq c \,.
  \label{primastimaa}
\Eeq

\step
Second a priori estimate

We take any $(v,\vG)\in\L2\calV$, test \rev{\eqref{primaaeps}} by $(v,\vG)$
and integrate over $(0,T)$.
By owing to \eqref{primastimaa}, we easily obtain~that
\Beq
  \intQ \dt(\peps+\tau\qeps) \, v
  + \intS \dt(\pGeps+\tau\qGeps) \, \vG
  \leq c \, \bigl( \norma v_{\L2V} + \norma\vG_{\L2\VG} \bigr) .
  \non
\Eeq
Equivalently, we have that
\Beq
  \ioT \< \dt(\peps+\tau\qeps,\pGeps+\tau\qGeps)(t) , (v,\vG)(t) >_\calV \, dt
  \leq c \, \norma{(v,\vG)}_{\L2\calV}
  \non
\Eeq
\anold{whereas it follows that}
\Beq
  \norma{\dt(\peps+\tau\qeps,\pGeps+\tau\qGeps)}_{\L2\calVp} \leq c \,.
  \label{secondastimaa}
\Eeq

\step
Conclusion of the existence proof

Thanks to \accorpa{primastimaa}{secondastimaa} and to standard compactness results,
we have~that 
\Bsist
  && (\peps,\pGeps) \to (p,\pG)
  \qquad \hbox{weakly star in $\L\infty\calV$}
  \non
  \\
  && (\qeps,\qGeps) \to \betti{(q,\qG)}
  \qquad \hbox{weakly star in $\L\infty\calH\cap\L2\calV$}
  \non
  \\
  && \dt(\peps+\tau\qeps,\pGeps+\tau\qGeps) \to \dt(p+\tau q,\pG+\tau\qG)
  \qquad \hbox{weakly in $\L2\calVp$}
  \non
\Esist
as $\eps\seto0$, 
for some quadruplet $(p,\pG,q,\qG)$ satisfying \Regsoluza, at least for a subsequence.
Moreover, \eqref{primastimaa} also implies~that
\Beq
  \eps \, \dt(\peps,\pGeps) \to (0,0)
  \qquad \hbox{strongly in $\L2\calH$} \,.
  \non
\Eeq
Then, it is \sfw\ to see that this limit quadruplet solves
equalities that are equivalent to \eqref{primaa} and \eqref{secondaa},
namely, their integrated versions with time dependent test functions $(v,\vG)\in\L2\calV$.
On the other hand, the above convergence properties
imply that $(\peps+\tau\qeps,\pGeps+\tau\qGeps)$ converges to $(p+\tau q,\pG+\tau\qG)$ weakly in $\C0\calH$.
Thus, the \anold{terminal} condition \eqref{cauchyaeps} and the \anold{\an{convergence} in} \eqref{hpzetaeps}
imply that \eqref{cauchya} is \anold{fulfilled} by $(p+\tau q,\pG+\tau\qG)$.
This completes the proof of the existence of a solution.

\step
Uniqueness

Since the problem is linear, it is sufficient to prove that the unique solution $\soluza$
to \Pbla\ is $(0,0,0,0)$
if all the functions $\zeta_i$ are replaced by zero\anold{s}.
\anold{In this direction, we } set for $t\in[0,T]$
\Bsist
  && P(t) := - \intT p(s) \, ds , \quad
  Q(t) := - \intT q(s) \, ds , \quad
  \Lam(t) := - \intT (\lam q)(s) \, ds
  \non
  \\
  && \PG(t) := - \intT \pG(s) \, ds , \quad
  \QG(t) := - \intT \qG(s) \, ds , \quad
  \LamG(t) := - \intT (\lamG\qG)(s) \, ds
  \non
\Esist
and \anold{integrate} \eqref{primaa} written at the time $s$ over $(t,T)$ with respect to~$s$.
Due to \eqref{cauchya}, it results that
\begin{align}
  & \< (p+\tau q,\pG+\tau\qG)(t) , (v,\vG) >_\calV
  - \iO \nabla Q(t) \cdot \nabla v
  - \iG \nablaG\QG(t) \cdot \nablaG\vG
  \non
  \\
  & = \gamma \iO P(t) v
  + \gamma \iG \PG(t) \vG
  + \iO \Lam(t) v 
  + \iG \LamG(t) \vG
  \label{intprimaa}
\end{align}
for every $t\in[0,T]$ and every $(v,\vG)\in\calV$.
Now, we test this equation by $(q,\qG)(t)$.
By avoiding writing the time $t$ for brevity, we obtain~that
\begin{align}
  & \iO pq
  + \tau \iO |q|^2
  + \iG \pG \qG
  + \tau \iG |\qG|^2
  - \frac 12 \, \frac d{dt} \iO |\nabla Q|^2
  - \frac 12 \, \frac d{dt} \iG |\nablaG\QG|^2
  \non
  \\
  & = \gamma \iO Pq
  + \gamma \iG \PG\qG
  + \iO \Lam q
  + \iG \LamG\qG \,.
  \label{test1}
\end{align}
Next, we test \eqref{secondaa} by $(p,\pG)$ and add the same quantity 
$-\frac12\,\frac d{dt}\iO|P|^2=-\iO Pp$
and the boundary analogue to both sides.
We get~that
\begin{align}
  & \iO |\nabla p|^2
  + \iG |\nablaG\pG|^2
  - \frac12 \,\frac d{dt} \iO |P|^2
  - \frac12 \,\frac d{dt} \iG |\PG|^2
  \non
  \\
  & = \iO qp
  + \iG \qG\pG
  - \iO Pp
  - \iG \PG\pG \,.
  \label{test2}
\end{align}
Finally, we test \eqref{intprimaa} by $\frac1{2\tau}(p,\pG)$ and have~that
\begin{align}
  & \frac 1{2\tau} \iO |p|^2
  + \frac 1{2\tau} \iG |\pG|^2
  \non
  \\
  & = - \frac 12 \iO qp
  - \frac 12 \iG \qG\pG
  + \frac 1{2\tau} \iO \nabla Q \cdot \nabla p
  + \frac 1{2\tau} \iG \nablaG\QG \cdot \nabla\pG
  \non
  \\
  & \quad {}
  + \frac \gamma{2\tau} \iO Pp
  + \frac \gamma{2\tau} \iG \PG\pG
  + \frac 1{2\tau} \iO \Lam p
  + \frac 1{2\tau} \iG \LamG\pG \,.
  \label{test3}
\end{align}
At this point, we add \accorpa{test1}{test3} to each other 
and notice that some cancellation occur.
Moreover, we integrate in time over $(t,T)$.
The sum of the bulk terms on the \lhs\ is given~by
\Beq
  S_l(t) := \tau \intQt |q|^2
  + \frac 12 \iO |\nabla Q(t)|^2 
  + \intQt |\nabla p|^2
  + \frac 12 \iO |P(t)|^2
  + \frac 1{2\tau} \intQt |p|^2 
  \label{lbulk}
\Eeq
while that on the \rhs\ is
\begin{align}
  & S_r(t) := \gamma \intQt Pq
  + \intQt \Lam q
  - \intQt Pp
  - \frac 12 \intQt qp
  \non
  \\
  & \quad {}
  + \frac 1{2\tau} \intQt \nabla Q \cdot \nabla p
  + \frac \gamma{2\tau} \intQt Pp
  + \frac 1{2\tau} \intQt \Lam p.
  \label{rbulk}
\end{align}
The most delicate term is the fourth integral in \eqref{rbulk},
for which sharp coefficients are needed in applying the Young inequality
(while for the other ones such \an{precision} is not necessary).
We thus have for every~$\delta>0$
(where the above terms are estimated in their order)
\Bsist
  && S_r(t)
  \leq \delta \intQt |q|^2
  + c_\delta \intQt |P|^2
  + \delta \intQt |q|^2
  + c_\delta \intQt |\Lam|^2
  \non
  \\
  && \quad {}
  + \delta \intQt |p|^2
  + c_\delta \intQt |P|^2
  + \frac \tau 4 \intQt |q|^2
  + \frac 1{4\tau} \intQt |p|^2
  \non
  \\
  && \quad {}
  + \frac 12 \intQt |\nabla p|^2
  + c \intQt |\nabla Q|^2
  + \delta \intQt |p|^2
  + c_\delta \intQt |P|^2
  + \delta \intQt |p|^2
  + c_\delta \intQt |\Lam|^2 \,.
  \non
\Esist
On the other hand, we have that
\Beq
  |\Lam(s)|^2 
  \leq \Bigl( \int_s^T c \, |q(s')| \, ds' \Bigr)^2
  \leq c \int_s^T |q(s')|^2 \, ds'
  \quad \hbox{for every $s\in[0,T]$}
  \non
\Eeq
whence also
\Beq
  \intQt |\Lam|^2
  \leq c \intT \Bigl( \iO \Bigl( \int_s^T |q(s')|^2 \, ds' \Bigr) \Bigr) \, ds
  = c \intT \Bigl( \int_{Q^s} |q|^2 \Bigr) \, ds \,.
  \non
\Eeq
Since the boundary terms coming from \accorpa{test1}{test3} are \anold{similar},
those of the \rhs\ can be estimated in the same way.
Therefore, by collecting everything, rearranging \anold{the terms}, choosing $\delta$ small enough, 
and applying the (backward) Gronwall lemma,
we conclude that $\soluza=(0,0,0,0)$.
\Edim

Finally, we express the necessary condition for optimality \eqref{preopt}
in terms of the solution to the adjoint problem.
As sketched in the Introduction, our result is the following:

\Bthm
\label{Opt}
Let $(\ustar,\uGstar)$ and $\solstar$ be an optimal control and the corresponding state, respectively,
and let $\soluza$ be the solution to the corresponding adjoint problem.
Then, the following variational inequality is satisfied
{%
\Beq
  \intQ (\gamma p + \alpha_5 \ustar) (u-\ustar)
  + \intS (\gamma\pG + \alpha_6 \uGstar) (\uG-\uGstar)
  \geq 0
  \label{opt}
\Eeq
for every $(u,\uG)\in\Uad$.}
\Ethm

\Bdim
We also consider the linearized system \Pbll\ associated \an{with} $\solstar$ and to the variation
$(h,\hG):=(u-\ustar,\uG-\uGstar)$.
We test \eqref{primal}, \eqref{secondal}, \eqref{primaa} and \eqref{secondaa} by
$(p,\pG)$, $(q,\qG)$, $-(\psi,\psiG)$ and $-(\eta,\etaG)$, respectively.
We obtain the following identities
\Bsist
  && \iO \dt\psi \, p
  + \iG \dt\psiG \, \pG
  + \iO \nabla\eta \cdot \nabla p
  + \iG \nablaG\etaG \cdot \nablaG\pG
  \non
  \\
  \separa
  && = \gamma\iO (u-\ustar-\psi) p
  + \gamma \iG (\uG-\uGstar-\psiG) \pG
  \non
  \\
  && \tau \iO \dt\psi \, q
  + \tau \iG \dt\psiG \, \qG
  + \iO \nabla\psi \cdot \nabla q
  + \iG \nablaG\psiG \cdot \nablaG\qG
  \non
  \\
  && \quad {}
  + \iO \lam \psi q
  + \iG \lamG \psiG \qG
  \non
  \\
  && = \iO \eta q
  + \iG \etaG \qG
  \non
  \\
  && \< \dt(p+\tau q,\pG+\tau\qG) , (\psi,\psiG) >_\calV
  \non
  \\
  && \quad {}
  - \iO \nabla q \cdot \nabla\psi
  - \iG \nablaG\qG \cdot \nablaG\psiG
  - \iO (\gamma p + \lam q) \psi
  - \iG (\gamma\pG + \lamG\qG) \psiG
  \non
  \\
  && = - \iO \zeta_1 \psi
  - \iG \zeta_2 \psiG
  \non
  \\
  && - \iO \nabla p \cdot \nabla\eta
  - \iG \nablaG\pG \cdot \nablaG\etaG
  \non
  \\
  && = - \iO q \eta
  - \iG \qG \etaG
  \non
\Esist
which hold \aet.
At this point, we sum up.
Almost all terms cancel each other and it remains~that
\Bsist
  && \< \dt(p+\tau q,\pG+\tau\qG) , (\psi,\psiG) >_\calV
  + \iO (p+\tau q) \dt\psi
  + \iG (\pG+\tau\qG) \dt\psiG
  \non
  \\
  && = \gamma\iO (u-\ustar) p
  + \gamma \iG (\uG-\uGstar) \pG
  - \iO \zeta_1 \psi
  - \iG \zeta_2 \psiG \,.
  \label{peropt}
\Esist
On the other hand, we recall that
\Beq
  (p+\tau q,\pG+\tau\qG) \in \H1\calVp \cap \L2\calV
  \aand
  (\psi,\psiG) \in \H1\calH \cap \L2\calV
  \non
\Eeq
whence, we can write
\Bsist
  && \< \dt(p+\tau q,\pG+\tau\qG) , (\psi,\psiG) >_\calV
  + \iO (p+\tau q) \dt\psi
  + \iG (\pG+\tau\qG) \dt\psiG
  \non
  \\
  && = \< \dt(p+\tau q,\pG+\tau\qG) , (\psi,\psiG) >_\calV
  + \bigl( (p+\tau q,\pG+\tau\qG) , \dt(\psi,\psiG) \bigr)_\calH
  \non
  \\
  && = \frac d{dt} \, \bigl( (p+\tau q,\pG+\tau\qG),(\psi,\psiG) \bigr)_\calH 
  \quad \aet.
  \non
\Esist
Hence, by integrating \eqref{peropt} over $(0,T)$ 
and accounting for the Cauchy conditions \eqref{cauchyl} and~\eqref{cauchya}, 
we deduce~that
\Bsist
  && \bigl( (\zeta_3,\zeta_4) , (\psi,\psiG)(T) \bigr)_\calH
  \non
  \\
  && = \gamma \intQ (u-\ustar) p
  + \gamma \intS (\uG-\uGstar) \pG
  - \intQ \zeta_1 \psi
  - \intS \zeta_2 \psiG \,.
  \non
\Esist
Equivalently, we have that
\Bsist
  && \intQ \zeta_1 \psi
  + \intS \zeta_2 \psiG 
  + \iO \zeta_3 \psi(T)
  + \iG \zeta_4 \psiG(T)
  \non
  \\
  && = \gamma\intQ (u-\ustar) p
  + \gamma \intS (\uG-\uGstar) \pG
  \non
\Esist
so that \eqref{opt} coincides with \eqref{preopt}.
\Edim


\section*{Acknowledgments}
\anold{
%
ER and AS gratefully acknowledge support 
from the MIUR-PRIN Grant 2020F3NCPX ``Mathematics for industry 4.0 (Math4I4)''
and their affiliation to the GNAMPA (Gruppo Nazionale per l'Analisi Matematica, 
la Probabilit\`a e le loro Applicazioni) of INdAM (Isti\-tuto 
Nazionale di Alta Matematica).
}


\vspace{3truemm}


\vspace{3truemm}

\Begin{thebibliography}{10}

\betti{%
\bibitem{BO}
M. Bahiana, Y. Oono, 
Cell dynamical system approach to block copolymers, 
{\em Phys. Rev. A} \textbf{41} (1990), 6763-6771.}

\bibitem{Barbu}
V. Barbu,
``Nonlinear Differential Equations of Monotone Type in Banach Spaces'',
Springer,
London, New York, 2010.

\betti{%
\bibitem{BEG}A.L. Bertozzi, S. Esedoglu, A. Gillette,
Inpainting of binary images using the Cahn--Hilliard equation,
{\em IEEE Trans. Image Process.} \textbf{16} (2007), 285-291.}

\betti{
\bibitem{BDS}
E. Bonetti, W. Dreyer, G. Schimperna, 
Global solution to a generalized Cahn--Hilliard
equation with viscosity, 
{\em Adv. Differential Equations} {\bf  8} (2003) 231-256.
}

\rev{
\bibitem{Brezis}
H. Brezis,
``Op\'erateurs maximaux monotones et semi-groupes de contractions
dans les espaces de Hilbert'',
North-Holland Math. Stud. Vol.
{\bf 5},
North-Holland,
Amsterdam,
1973.}

\betti{
\bibitem{CaCo} 
L. Calatroni, P. Colli,
Global solution to the Allen--Cahn equation with singular potentials and
dynamic boundary conditions, 
{\em Nonlinear Anal.} {\bf 79} (2013), 12-27.}

\betti{%
\bibitem{CRW}
C. Cavaterra, E. Rocca, H. Wu, 
Long-time dynamics and optimal control of a diffuse interface model for tumor growth,
{\em Appl. Math. Optim.} {\bf 83} (2021), 739-787.}

\betti{
\bibitem{CF}
P. Colli, T. Fukao, 
Cahn–Hilliard equation with dynamic boundary conditions and mass constraint on the boundary, 
{\em J. Math. Anal. Appl.} {\bf 429} (2015), 1190-1213.
}

\betti{
\bibitem{CGH}
P. Colli, G. Gilardi, D. Hilhorst,
On a Cahn--Hilliard type phase field system related to tumor growth,
{\em Discret. Cont. Dyn. Syst.} {\bf 35} (2015), 2423-2442.}


\betti{%
\bibitem{CGRS}
P. Colli, G. Gilardi, E. Rocca, J. Sprekels, 
Optimal distributed control of a diffuse interface model of tumor growth,
{\em Nonlinearity}  {\bf 30} (2017), 2518-2546.}

\bibitem{CGRS4}
P. Colli, G. Gilardi, E. Rocca, J. Sprekels,
Well-posedness and optimal control for a Cahn--Hilliard--Oono system
with control in the mass term,
\betti{ {\em Discret. Contin. Dyn. Syst. Ser. S}  {\bf 15} (2022), 2135-2172.}

\betti{
\bibitem{CGS2014}
\anold{P. Colli, G. Gilardi, E. Rocca,}
 On the Cahn–Hilliard equation with dynamic boundary conditions and a dominating boundary potential, 
{\em J. Math. Anal. Appl.} {\bf 419} (2014), 972-994.}

\betti{
\bibitem{CGS2017}
P. Colli, G. Gilardi, J. Sprekels, 
Recent results on the Cahn--Hilliard equation with dynamic boundary conditions, 
{\em Vestn. Yuzhno-Ural. Gos. Univ., Ser. Mat. Model. Program.} {\bf 10} (2017), 5-21.}

\bibitem{CGS13}
P. Colli, G. Gilardi, J. Sprekels,
On a Cahn--Hilliard system with convection and dynamic boundary conditions,
{\it  Ann. Mat. Pura Appl.} {\bf 197} (2018), 1445-1475. 

\bibitem{CGS14}
P. Colli, G. Gilardi, J. Sprekels,
Optimal velocity control of a viscous Cahn--Hilliard system 
with convection and dynamic boundary conditions,
{\it SIAM J. Control Optim.} {\bf 56} (2018), 1665-1691. 


\betti{\bibitem{CGS-2019b}
P. Colli, G. Gilardi, J. Sprekels,
Optimal velocity control of a convective Cahn--Hilliard system with double obstacles and 
dynamic boundary conditions:\ a `deep quench' approach,
{\em J. Convex Anal.} {\bf 26} (2019), 485-514.}
%


\betti{
\bibitem{CS}
P. Colli, A. Signori,
Boundary control problem and optimality conditions for the Cahn--Hilliard equation with dynamic boundary conditions,
{\em Internat. J. Control} {\bf  94} (2021), 1852–1869. 
}

\anold{
\bibitem{CSS2}
P. Colli, A. Signori, J. Sprekels,
Second-order analysis of an optimal control problem in a phase field tumor growth model with singular potentials and chemotaxis,
{\it ESAIM Control Optim. Calc. Var.}, {\bf 27} (2021). }

\betti{
\bibitem{ES}
C.M. Elliott, A.M. Stuart, 
Viscous Cahn--Hilliard equation. II. Analysis, 
{\em J. Differential
Equations} {\bf 128} (1996), 387-414.}

\betti{
\bibitem{FMD1}
H.P. Fischer, Ph. Maass, W. Dieterich, 
Novel surface modes in spinodal decomposition,
{\em Phys. Rev. Letters} {\bf 79} (1997), 893-896.}

\betti{
\bibitem{FMD2}
H.P. Fischer, Ph. Maass, W. Dieterich, 
Diverging time and length scales of spinodal
decomposition modes in thin flows, 
{\em Europhys. Letters} {\bf 42} (1998), 49-54.}

\betti{%
\bibitem{FGR}
S.~Frigeri, M.~Grasselli, E.~Rocca,
On a diffuse interface model of tumor growth,
{\em European J. Appl. Math.} {\bf 26} (2015), 215-243.}

\rev{
\bibitem{FN}
T. Fukao, Y. Noriaki,
A boundary control problem for the equation and dynamic boundary condition of Cahn–Hilliard type,
{\it Solvability, Regularity, and Optimal Control of Boundary Value Problems for PDEs: In Honour of Prof. Gianni Gilardi}
(2017), 255-280.
}

\betti{
\bibitem{GK}
H. Garcke, P. Knopf, 
Weak solutions of the Cahn–Hilliard system with dynamic boundary conditions: a gradient flow approach,
{\em SIAM J. Math. Anal.} {\bf 52} (2020). }

\bibitem{garcke}
\anold{
H.~Garcke, K.~F. Lam, E.~Sitka, V.~Styles.
\newblock A {C}ahn--{H}illiard--{D}arcy model for tumour growth with chemotaxis
  and active transport.
\newblock {\em Math. Models Methods in Appl. Sci.}, {\bf 26}(06) (2016), 1095-1148.}

\betti{\bibitem{GARL_1}
H. Garcke, K. F. Lam,
Well-posedness of a Cahn--Hilliard system modelling tumour
growth with chemotaxis and active transport,
{\em European. J. Appl. Math.} {\bf 28} (2017), 284--316.}

\betti{%
\bibitem{GLR}
H. Garcke, K. F. Lam, E. Rocca,
Optimal control of treatment time in a diffuse interface model of tumor growth,
{\em  Appl. Math. Optim.} {\bf 78} (2018), 495-544.}

\bibitem{GiMiSchi} 
G. Gilardi, A. Miranville, G. Schimperna,
On the Cahn--Hilliard equation with irregular potentials and dynamic boundary conditions,
{\em  Commun. Pure Appl. Anal.\/} 
{\bf 8} (2009), 881-912.

\bibitem{GiMiSchi2} 
G. Gilardi, A. Miranville, G. Schimperna,
Long-time behavior of the Cahn-Hilliard equation with irregular potentials and dynamic boundary conditions,
{\em Chinese Annals of Mathematics Series~B\/} {\bf 31} (2010) 679-712. 

\betti{
\bibitem{GGM}
A. Giorgini, M. Grasselli, A. Miranville,
The Cahn--Hilliard--Oono equation with singular potential,
{\em Math. Models Methods Appl. Sci.}
{\bf 27} (2017), 2485-2510.}

\betti{
\bibitem{Ketal} 
R. Kenzler, F. Eurich, Ph. Maass, B. Rinn, J. Schropp, E. Bohl, W. Dieterich, 
Phase
separation in confined geometries: solving the Cahn--Hilliard equation with generic boundary
conditions, 
{\em Comput. Phys. Comm.} {\bf 133} (2001), 139-157.}

\betti{%
\bibitem{KS}
E. Khain, L.M. Sander,
Generalized Cahn--Hilliard equation for biological applications, 
{\em Phys. Rev. E} \textbf{77} (2008), 051129-1-051129-7.}

\anold{
\bibitem{KS2}
P. Knopf, A. Signori, Existence of weak solutions to multiphase Cahn--Hilliard--Darcy and Cahn--Hilliard--Brinkman models for stratified tumor growth with chemotaxis and general source terms,
{\it Comm. Partial Differential Equations}, {\bf 47}(2) (2022), 233-278.
}

\anold{
\bibitem{KSnon}
P.~Knopf, A.~Signori,
On the nonlocal Cahn--Hilliard equation with nonlocal dynamic boundary condition and boundary penalization,
{\it J. Differential Equations}, {\bf 280}(4) (2021), 236-291. 
}

\bibitem{Lions}
J.-L.~Lions,
``Quelques m\'ethodes de r\'esolution des probl\`emes
aux limites non lin\'eaires'',
Dunod; Gauthier-Villars, Paris, 1969.

\betti{
\bibitem{LW}
C. Liu, H. Wu, An energetic variational approach for the Cahn--Hilliard equation with dynamic boundary conditions: model derivation and mathematical analysis, 
{\em Arch. Ration. Mech. Anal.} {\bf  233 }(2019), 167-247.}

\betti{
\bibitem{MMP}
A. Makki, A. Miranville, M. Petcu, 
A numerical analysis of the coupled Cahn--Hilliard/Allen--Cahn system with dynamic boundary conditions.,
{\em Int. J. Numer. Anal. Model.} {\bf  19} (2022), 630–655.} 

\betti{%
\bibitem{MR}
S. Melchionna, E. Rocca,
On a nonlocal Cahn--Hilliard equation with a reaction term, 
{\em Adv. Math. Sci. Appl.} {\bf 24} (2014), 461-497.}

\betti{%
\bibitem{M}
A. Miranville, Asymptotic behavior of the Cahn--Hilliard--Oono equation, 
{\em J. Appl. Anal. Comput.} {\bf 1} (2011), 523-536.}

\betti{%
\bibitem{MRS} 
A. Miranville, E. Rocca, G. Schimperna, 
On the long time behavior of a tumor growth model, 
{\it J. Differential Equations} {\bf 67} (2019),  2616-2642.}

\bibitem{MiZe}
A. Miranville, S. Zelik, 
Robust exponential attractors for Cahn--Hilliard type equations with singular potentials, 
{\em Math. Methods Appl. Sci.} {\bf 27} (2004), 545-582.

\betti{
\bibitem{N}
A. Novick-Cohen, 
Energy methods for the Cahn--Hilliard equation, 
{\em Quart. Appl. Math.} {\bf  46}
(1988), 681-690.}

\betti{%
\bibitem{oonopuri87}
Y. Oono, S. Puri, 
Computationally efficient modeling of ordering of quenched phases, 
{\em Phys. Rev. Lett.} {\bf 58} (1987), 836-839.}

\betti{%
\bibitem{oonopuri88I}
Y. Oono, S. Puri, 
Study of phase-separation dynamics by use of cell dynamical systems. I. Modeling,
{\em Phys. Rev. A} {\bf 38} (1988), 434-453.}

\betti{%
\bibitem{oonopuri88II}
S. Puri, Y. Oono,
Study of phase-separation dynamics by use of cell dynamical systems. II. Two-dimensional demonstrations, 
{\em Phys. Rev. A} {\bf 38} (1988), 1542-1565.}

%

\anold{
\bibitem{RSchS} 
E.~Rocca, G.~Schimperna, A.~Signori,
On a Cahn--Hilliard--Keller--Segel model with generalized logistic source describing tumor growth,
{\it J. Differential Equations}, {\bf 343} (2023),  530-578.
}

\anold{
\bibitem{SS}
L.~Scarpa, A.~Signori,
On a class of non-local phase-field models for tumor growth with possibly singular potentials, chemotaxis, and active transport,
{\it Nonlinearity}  {\bf 34} (2021), 3199-3250. 
}

\anold{
\bibitem{S}
A. Signori,
Optimal distributed control of an extended model of tumor
growth with logarithmic potential.
{\it Appl. Math. Optim.} {\bf 82} (2020), 517-549.
}

\End{thebibliography}

\End{document}
